\documentclass[11pt,leqno]{article}
\usepackage{graphicx, amsfonts, amsthm, amsxtra, amssymb, verbatim,rotating}
\usepackage[mathscr]{euscript}
\def\GL{{\rm GL}}
\def\ZZ{{\Bbb Z}}

\begin{document}

\begin{center}
{\LARGE\bf Halphen's transform and middle convolution}

\vspace{.25in} 
{\large {\sc Stefan Reiter}}
\end{center}

\begin{abstract}
We show that the Halphen transform of a Lam\'e equation can be written as the symmetric square of the Lam\'e equation followed by an Euler transform. 
We use  this to compute a list of Lam\'e equations with arithmetic Fuchsian monodromy group.  It contains all those Lam\'e equations
where the  quaternion algebra $A$ over $k$ associated 
to the arithmetic Fuchsian group
is a  quaternion algebra $A$ over $\mathbb{Q}$.
Further we classify all geometric braid group orbits in $\rm{SL}_2(\ZZ)^4$
with the possible exception of three orbits.
\end{abstract}

\def\Gal{{\rm Gal}}              
\def\ZZ{\mathbb{Z}}                   
\def\QQ{\mathbb{Q}}                   
\def\CC{\mathbb{C}}                   
\def\ring{{\sf R}}                             
\def\RR{\mathbb{R}}                   
\def\NN{\mathbb{N}}                   
\def\uhp{{\mathbb H}}                
\newcommand{\mat}[4]{
     \begin{pmatrix}
            #1 & #2 \\
            #3 & #4
       \end{pmatrix}
    }                                

\def\per{{\sf pm}}

\def\la{\lambda}
\def\th{\theta}
\def\P{\mathbb P}
\theoremstyle{plain}
\def\PP{{\mathbb P}}
\def\SL{{\rm SL}}
\def\SU{{\rm SU}}
\def\GO{{\rm GO}}
\def\SO{{\rm SO}}
\def\Sp{{\rm Sp}}
\def\GL{{\rm GL}}
\def\PGL{{\rm PGL}}
\def\dz{{\rm d} z}
\def\dx{{\rm d} x}
\def\A{{\bf A}}
\def\B{{\bf B}}
\def\C{{\bf C}}
\def\T{{\bf T}}
\def\la{{\lambda}}
\def\nlambda{{\tilde{\lambda}}}
\def\nue{{\nu}}
\def\ele{{\rm L}}
\def\tparam{{t}}
\def\Mu{{\mu}}
\def\th{{\theta}}
\def\thh{{\tilde{\theta}}}
\def\Th{{\Theta}}
\def\tr{{\rm tr}}
\def\rk{{\rm rk}}
\def\Trace{{\rm Tr}}
\def\Mat{{\rm Mat}}
\def\al{{\alpha}}
\def\diag{{\rm diag}}
\def\id{{\rm id}}
\def\k{{\frak k}}
\def\l{{\frak  l}}
\def\Res{{\rm  Res}}
\def\trace{{\rm  trace}}
\def\Id{{\rm Id}}
\def\ten{{\otimes}}
\def\si{{\sigma}}
\def\sym{{\rm sym }}
\def\antidiag{{\rm antidiag }}

\renewcommand{\mat}[4]{
     \begin{pmatrix}
            #1 & #2 \\
            #3 & #4
       \end{pmatrix}
    }                               
\newcommand{\matv}[1]{
     \begin{pmatrix}#1
      \end{pmatrix}}

{

\newtheorem{thm}{Theorem.}[section]
\newtheorem{cor}[thm]{Corollary.}
\newtheorem{lemma}[thm]{Lemma.}
\newtheorem{prop}[thm]{Proposition.}
\newtheorem{introthm}{Theorem}}

{

\theoremstyle{definition}
\newtheorem{defn}[thm]{Definition.}
\newtheorem{prob}[thm]{Problem}
\newtheorem{rem}[thm]{Remark.}
\newtheorem{ex}[thm]{Example}
\newtheorem{ass}[thm]{Assumption}
\newtheorem{notation}[thm]{Notation.}
\newtheorem{claim}[thm]{Claim}
\newtheorem{cond}[thm]{Condition}
\newtheorem{conj}[thm]{Conjecture}
\newtheorem{question}[thm]{Question}
\newtheorem{sit}[thm]{Situation}
\newtheorem{var}[thm]{Variant}
\newtheorem{constr}[thm]{Construction}
\newtheorem{hyp}[thm]{Hypothesis}
\newtheorem{tab}[thm]{Table}

}

\renewcommand{\em}{\sl}
\newcommand{\eps}{\epsilon}
\makeatletter
\renewcommand{\subsection}{\@startsection{subsection}{2}%
{\z@}{-3.25ex plus -1ex minus-.2ex}{-1em}{\bf}} \makeatother

\renewcommand{\labelenumi}{{\rm (\roman{enumi})}}
\parindent0em
\parskip0.2ex

\renewcommand{\labelenumi}{{\rm (\roman{enumi})}}

\newcommand{\textcyr}[1]{{\fontencoding{OT2}\selectfont #1}}


\section{Introduction}

Besides the Gauss hypergeometric differential equation, the Lam\'e equation
\begin{eqnarray*}   p(x)y''+\frac{1}{2}p'(x) y'+q(x)y&=&0,  \\
p(x) = 4\prod_{i=1}^3 (x-e_i)=4 x^3-g_2x-g_3,&
 q(x)=&- ( n(n+1)x-H) \end{eqnarray*}
 is one of the best studied  second order differential equation.
Of special interest are those  Lam\'e equations with finite
monodromy group, having therefore algebraic solutions, studied by
Baldassarri, Beukers and van der Waall, Chudnovsky and Chudnovsky, Dwork and
many others (cf. \cite{Ba}, \cite{BW}, \cite{Chud1383} and \cite{Dwork90}).

Lam\'e equations also occur in the context of  Grothendieck's $p$-curvature conjecture (cf. \cite[p. 15]{Chud1383}). 
This conjecture
says that if  the $p$-curvature of a differential equation
is zero modulo $p$ for almost all primes $p$
then its monodromy group is finite. More generally, it is conjectured that
if the  $p$-curvature is globally 
nilpotent then the differential equation is { \it geometric } (also called coming from geometry, i.e. it is a  product of factors of  Picard-Fuchs differential equations, see \cite[Chap. II \S 1]{Andre89}).

These conjectures are proven by Chudnovsky and Chudnovsky
in the Lam\'e case for $n$ being an integer. In this case
the monodromy group is a dihedral group or reducible
(s. \cite[Thm. 2.1]{Chud1383}).
Moreover they showed that, for a given Riemann scheme, there is only a finite number of Lam\'e operators that are globally nilpotent (s. \cite[Thm. 2.3]{Chud1383}). (By a result
of Deligne in \cite{Del} there are only finitely many geometric differential equations with given Riemann scheme.)

One also knows that  Lam\'e equations  
with an arithmetic Fuchsian monodromy
group of signature $(1,e)$ provide interesting examples of
geometry differential equations (s. \cite{Chud1383}). 
In \cite{Krammer96} Krammer  determined one such example and showed
that it is not a (weak) pull-back of a hypergeometric differential equation
contradicting a conjecture of Dwork that any globally nilpotent second order differential equation on $\PP^1/\bar{\QQ}$ has either algebraic solutions, or is a weak pull-back of a Gauss hypergeometric differential equation  (cf. \cite[Section 11]{Krammer96}).
  Recently  Sijsling classified in \cite{Sijsling} all  Lam\'e equations  with arithmetic Fuchsian monodromy group of signature $(1,e)$ that are pull-backs of hypergeometric
differential equations.

But there is also the Halphen transform that changes the Lam\'e equation
into another second order differential equation, again a Heun equation.
This was used in \cite{Chud1383} for $n=-\frac{1}{2}.$
In this case it turned out that the new differential equation
is a pull-back of  a hypergeometric differential equation.

In this paper we will make use of the last observation.
We show in Section~\ref{Halphen}, Corollary~\ref{Mu}, 
that the  Halphen transform of a Lam\'e equation can be written as the symmetric square of the Lam\'e equation followed by an Euler transform.
Since these are geometric operations (cf. \cite{Andre89}) 
the obtained Heun equation
is also geometric provided the Lam\'e equation is it. 
We also generalize the Halphen transform to the case where 
$n$ is not necessarily $-\frac{1}{2}$.

 In Section~\ref{Ex} we go the opposite way, by starting 
with  special geometric Heun equations, that were
jointly with H. Movasati computed in  \cite{MR08-2}, listed in { Table~3.1}:

\[
\begin{array}{llll}
   \mbox{Nr. in} &  p(x)y''+p'(x)y'+(\alpha \beta x+\tilde{H})y=0& p(x)\\
  \cite{MR08-2}   &\\
  1 & p(x)y''+p'(x)y'+(x+1/3)y=0 &x(x^2+x+1/3) \\
 2 &  p(x)y''+p'(x)y'+xy=0 &  x(x-1)(x+1) \\
 3 &  p(x)y''+p'(x)y'+(x-1/4)y=0& x(x-1)(x+1/8)\\
  4&  p(x)y''+p'(x)y'+(x+3)y =0 &x(x^2+11x-1) \\
  7&  p(x)y''+p'(x)y'+(35/36x - 3/2)y =0 &  x(x^2-14/3x+9) \\
 8 & p(x)y''+p'(x)y'+(35/36x-9/8)y=0&   x(x-1)(x-81/32) \\
 9 &  p(x)y''+p'(x)y'+(35/36x-45/2)y=0,&   x(x-1)(x-81) \\
 10&  p(x)y''+p'(x)y'+(35/36x-20/27)y=0,&   x(x-1)(x-32/27) \\
  11&  p(x)y''+p'(x)y'+(15/16x +1/8)y=0 & x(x^2+13/32x+1/8)\\
12&   p(x)y''+p'(x)y'+(15/16 x-3/2)y=0,&   x(x-1)(x-4) \\
 13&   p(x)y''+p'(x)y'+(15/16 x-12)y=0,&   x(x-1)(x-128/3) \\
 15 &   p(x)y''+p'(x)y'+(8/9 x)y=0,&   x(x-1)(x+1) \\
 16 &  p(x)y''+p'(x)y'+(8/9x-8/27)y=0,&   x(x-1)(x-2/27) \\
\end{array}
\]

Thus we obtain the { Table~3.2} of geometric Lam\'e equations  with
arithmetic Fuchsian monodromy group:

\[
\begin{array}{lllll}
   \mbox{Nr.} &
    p(x)y''+\frac{1}{2}p'(x)y'- ( n(n+1)x-H)y=0& p(x)=\\
 1  & p(x)y''+\frac{1}{2}p'(x)y'+(1/4x+1/12)y=0 &4x(x^2+x+1/3) \\
 2 &  p(x)y''+\frac{1}{2}p'(x)y'+ 1/4xy=0 &  4x(x-1)(x+1) \\
 3  &  p(x)y''+\frac{1}{2}p'(x)y'+(1/4 x-1/32)y=0& 4x(x-1)(x+1/8)\\
  4  & p(x)y''+\frac{1}{2}p'(x)y'+(1/4x+1/4)y=0  &4x(x^2+11x-1) \\
 7     & p(x)y''+\frac{1}{2}p'(x)y'+(2/9x - 1/3)y =0 &  4x(x^2-14/3x+9) \\
 8 &   p(x)y''+\frac{1}{2}p'(x)y'+(2/9x-31/128)y=0&  4 x(x-1)(x-81/32) \\
 9 &  p(x)y''+\frac{1}{2}p'(x)y'+(2/9x-2)y=0,&  4 x(x-1)(x-81) \\
 10&      p(x)y''+\frac{1}{2}p'(x)y'+(2/9x-7/36)y=0,& 4  x(x-1)(x-32/27) \\
 11&    p(x)y''+\frac{1}{2}p'(x)y'+(3/16x +3/128)y=0 & 4x(x^2+13/32x+1/8)\\
 12 &   p(x)y''+\frac{1}{2}p'(x)y'+(3/16 x-1/4)y=0,&  4 x(x-1)(x-4) \\
 13  &   p(x)y''+\frac{1}{2}p'(x)y'+(3/16 x-13/12)y=0,&  4 x(x-1)(x-128/3) \\
 15 &    p(x)y''+\frac{1}{2}p'(x)y'+(5/36 x)y=0,&  4 x(x-1)(x+1) \\
 16 &   p(x)y''+\frac{1}{2}p'(x)y'+(5/36x-1/36)y=0,& 4  x(x-1)(x-2/27) 
\end{array}
\]

Here we have considered the case when the
monodromy group of the above Heun equations 
is contained in $\SL_2(\ZZ)$ and has 
at least $3$ unipotent monodromy group generators.
These equations
arise from the  classification of families of elliptic curves
having $4$ singular fibres in \cite{Herfurtner}. 
Thus being rational
pull-backs of a (geometric) Gauss hypergeometric differential equation 
they are geometric.

In general
the Euler transform does not commute with pull-backs and in general destroys
properties of the monodromy group like being arithmetic or even being discrete
(s. e.g.  \cite{DR06-1}).
Hence most of the  
 Lam\'e equations obtained in the above way are not pull-backs of 
a Gauss hypergeometric differential equation.

 In the literature we found only  
  Lam\'e equations with
 arithmetic Fuchsian monodromy group (s. \cite{Krammer96} and \cite{Chud1383})
 or geometric Heun equations with 
 nondiscrete  monodromy group (s. \cite{DR06-1})
 being no pull-backs of hypergeometric
 differential equations, providing counter examples to the above mentioned  conjecture of Dwork.
  (Recently also  differential equations with $5$ singularities and
 non-arithmetic Fuchsian monodromy group were computed in \cite{BoMo}.)
 The Lam\'e equations having an arithmetic Fuchsian monodromy group 
 with signature $(1,e)$
 were computed by Chudnovsky and Chudnovsky
 but they listed only a few interesting cases
 in \cite{Chud1383}. 
 In addition they are obtained 
 via numerical solutions of the uniformization problem of punctured tori.
 However via our approach   one could perhaps modify Dwork's conjecture  in the
following way:
 Any globally nilpotent second order differential equation on $\PP^1/\bar{\QQ}$ is related to  a Gauss hypergeometric differential equation  via geometric operations
 (Beukers  related Krammer's example to a Lauricella  hypergeometric function of type FD and saved so Dwork's conjecture).

 In Section~\ref{Mono} we determine the monodromy group generators (up to braid group action)
of the Heun equations in Table~3.1  and show how to obtain the corresponding
  monodromy group generators of
 the Lam\'e equations in Table~3.2. 
 This can be done by classifying minimal tuples
 in the braid group orbit via the Fricke relation.

In  Section~\ref{Hodge} we derive some properties of
monodromy groups of geometric differential equations
from Deligne's paper \cite{Del}.
Using operations that preserve geometric differential equations,
like rational pull-backs, tensor products and the middle convolution
we obtain criteria to rule out that certain monodromy group generators
arise from geometric differential equations.
As an application we reprove the already mentioned result of
Chudnovsky and Chudnovsky  for the $p$-curvature conjecture
for Lam\'e equations with integer $n$ (see \cite[Thm. 2.1]{Chud1383}).

 In the appendix  we classify
 all braid group orbits of four monodromy group generators
 in $\SL_2(\CC)$ with integer Fricke parameters.
 We call a  braid group orbit {\it geometric} if it contains
 a tuple of monodromy group generators of a Picard-Fuchs
equation.
 A geometric monodromy group preserves a
 hermitian form, a result due to Deligne, see \cite{Del}. 
 We are indebted to Duco van Straten
 for communicating this crucial reference. Further the middle convolution
 plays an important role since it allows to construct from monodromy group
 generators of second order differential equations
  monodromy group
 generators of
 differential equations of higher order with orthogonal monodromy group.
 Using exceptional isomorphisms of Lie groups of small rank 
 \[ \sym^2 \SL_2(\CC) = \SO_3(\CC),\quad \Lambda^2 \Sp_4(\CC)= \SO_5(\CC),\quad
     \Lambda^2 \SL_4(\CC) = \SO_6(\CC)\]
 we obtain restrictions
 for the signature of the hermitian form and restrictions
 for the traces of certain elements in the monodromy group in the case
of a geometric braid group orbit.
 This allows in our case
 to distinguish the non geometric and geometric braid group orbits
 in $\SL_2(\ZZ)$ up to three open cases.
 It turns out that all these  geometric braid group orbits contain
 a monodromy tuple of a geometric Heun equation that 
 arises from the classification
 of families of elliptic curves having four singular fibres \cite{Herfurtner}.

\section{The Halphen transform of a Lam\'e equation}\label{Halphen}

We recall some properties of the Halphen transform
taken from  \cite[chap IX]{Poole} and \cite[p. 60-62]{Chud1240}.
The Lam\'e equation can be written in the algebraic form
\begin{eqnarray*}   p(x)y''+\frac{1}{2}p'(x) y'+q(x)y&=&0,  \\
p(x) = 4\prod_{i=1}^3 (x-e_i)=4 x^3-g_2x-g_3,&
 q(x)=&- ( n(n+1)x-H) \end{eqnarray*}
with  Riemann scheme
\begin{eqnarray*}
 \left(\begin{array}{cccc}
    e_1 &  e_2 &  e_3 & \infty \\
     0 &  0  & 0 & -\frac{n}{2}  \\
    \frac{1}{2}&  \frac{1}{2}&  \frac{1}{2}&  \frac{n+1}{2}
 \end{array}\right)
\end{eqnarray*}
or in the elliptic form
\begin{eqnarray*} y''-(n(n+1){\frak{p}(u)}-H)y&=&0,\end{eqnarray*}
where ${\frak{p}}(u)$ denotes the Weierstrass ${\frak{p}}(u)$-function.
The Halphen transform is obtained   by putting  
\begin{eqnarray*} u=2v,& y={\frak{p}'(v)}^{-n}z. \end{eqnarray*}
Then $z$ satisfies
\begin{eqnarray*} z''-2n\frac{{\frak{p}}''(v)}{{\frak{p}}'(v)}z'+
        4 (n (2n-1)\frak{p}(v) +H)z&=&0 . \end{eqnarray*}
Putting $x={\frak{p}}(v)$  we get again the algebraic form
\begin{eqnarray*}  p(x)y''+(\frac{1}{2}-n)p'(x)y'+4(n(2n-1)x+H)y&=&0\end{eqnarray*}
with Riemann scheme
\begin{eqnarray*}
 \left(\begin{array}{cccc}
    e_1 &  e_2 &  e_3 & \infty \\
     0 &  0  & 0 & -2n  \\
     n+\frac{1}{2}&  n+\frac{1}{2}&  n+\frac{1}{2}&  \frac{1}{2}-n
 \end{array}\right).
\end{eqnarray*}

We will show that the Halphen transform 
of a Lam\'e equation can be written as the symmetric square of the Lam\'e equation followed by an Euler transform.

Thus we recall the 

\begin{rem} (\cite[p. 59]{Chud1240})
The
symmetric square of a second order differential equation
 \begin{eqnarray} \label{sym0}  y''+q_1(x)y'+q_2(x)y&=&0  \end{eqnarray}
can be written as
\begin{eqnarray} \label{sym1}  y'''+3 q_1(x) y'' +(q_1'(x)+4 q_2(x) +2 q_1(x)^2) y'+2(q_2'(x)+2q_1(x)q_2(x))y&=&0. \end{eqnarray}
(I.e. all products of solutions of \eqref{sym0} satisfy  \eqref{sym1}).
Moreover in the Lam\'e case we have 
\begin{eqnarray*}p(x)=4(x-e_1)(x-e_2)(x-e_3),\quad q_1(x)=\frac{p'(x)}{2p(x)},&&q_2(x)=\frac{- (n(n+1)x-H)}{p(x)}\end{eqnarray*}
and therefore the symmetric square of a Lam\'e equation is 
\begin{eqnarray}\label{symsquare}
  p(x) y'''+\frac{3}{2} p'(x) y'' +(\frac{p''(x)}{2}- 4( n(n+1)x-H)) y'-
2 n(n+1)y&=&0.\end{eqnarray}
\end{rem}

Using the formula for the Euler integral in \cite[Chap. 3.3, 3.4]{IKSW} we get
the following

\begin{lemma}
Let $f$ be a solution of (\ref{symsquare}).
Then the Euler integral $\int_{\gamma}f(t) (x-t)^{-1-\mu} dt$
over a Pochhammer
double loop $\gamma$
satisfies
\begin{eqnarray*}  p(x)y'''+(\frac{3}{2}+\mu) p'(x)y''+r_1(x)y'+r_0(x)y&=&0, \end{eqnarray*} 
where 
\begin{eqnarray*}
p(x)&=&4(x-e_1)(x-e_2)(x-e_3),\\
    r_1(x)&=& 4 (6 x +2(-\sum e_i))( \frac{\mu (\mu-1)}{2}+\frac{3}{2} \mu +\frac{1}{2})
    -4(n(n+1)x-H), \\
 r_0(x)&=& 4 \mu^3+6\mu^2+2\mu-4 \mu (n(n+1))-2 n(n+1)=2(2\mu+1)(\mu-n)(\mu+n+1).
\end{eqnarray*}
\end{lemma}

Thus if we choose $\mu$ such that $r_0(x)=0$ then we get again a second order
differential equation: 

\begin{cor}\label{Mu}
We get the following special cases in the above lemma
\begin{enumerate}
\item[a)] If $\mu=-\frac{1}{2}$ then
\begin{eqnarray*} p(x) y''+ p'(x)y'+  (4H-\sum_{i=1}^3 e_i-(2n+3)(2n-1)x)y&=&0. \end{eqnarray*}

\item[b)] If $\mu=n$ then
\begin{eqnarray*} p(x)y''+ (\frac{3}{2}+n) p'(x)y'+(4H-(n+1)^2\sum_{i=1}^3 e_i+4(n+1)(2n+3)x)y&=&0.\end{eqnarray*}
 
\item[c)] If $\mu=-n-1$ then
\begin{eqnarray*} p(x)y''+ (\frac{1}{2}-n) p'(x)y'+(4(H-n^2\sum_{i=1}^3 e_i)+4n(2n-1)x)y&=&0. \end{eqnarray*}

\end{enumerate}
The case c) gives the Halphen transform, since we had there assumed $\sum_{i=1}^3 e_i=0$.
\end{cor}

Since convolution and tensor products are geometric operations
(cf. \cite[Chap. II \S 1]{Andre89}) we get the following

\begin{cor}
The Halphen transform preserves geometric differential equations.
\end{cor}

\begin{rem}
Applying the Euler integral and factoring out trivial subspaces is
exactly  the middle convolution operation (s. \cite{DR07}).
\end{rem}

We apply Corollary~\ref{Mu} a) to the following example studied 
by Krammer in \cite{Krammer96}, which is considered
as a counterexample to a conjecture of Dwork, 
being not a (weak) pull-back of a Gauss hypergeometric differential equation: 
(It also appears in \cite[p. 23]{Chud1383})

\begin{ex}
The geometric Lam\'e equation
 with arithmetic monodromy group of signature $(1,3)$
\[p(x)  y''+ \frac{1}{2}p'(x)y'+(\frac{2}{9} x-2)y=0, \quad p(x)=4 x(x-1)(x-81)\]
becomes after applying Corollary~\ref{Mu} a)
the Heun equation
\[ p(x) y''+p'(x)y'+(\frac{35}{9}x-90)y=0\]
with unipotent local monodromy at $0,1,$ and $81$.
We will see in the next section (Table~3.1, row 9) that this differential equation
is a rational pull-back of a hypergeometric differential equation.
\end{ex}

\begin{rem}
 The Halphen transform of a Lam\'e equation can also be written  as a rational pull-back
 $\phi(x)$ followed by a multiplication with an algebraic function. 
 E.g. if $(e_1,e_2,e_3)=(0,1,t)$  then the pull-back is
 \[ \phi(x)=\frac{(x^2-t)^2}{4x(x-1)(x-t)}\]
 and the algebraic function is
 \[ (4x(x-1)(x-t))^{n/2}. \]
\end{rem}

\begin{proof}
 Let $f(x)$ be a solution of the Lam\'e equation
 \begin{eqnarray*}   p(x)y''+\frac{1}{2}p'(x) y'+q(x)y&=&0,  \\
p(x) = 4\prod_{i=1}^3 (x-e_i)=4 x(x-1)(x-t),&
 &q(x)=- ( n(n+1)x-H). \end{eqnarray*}
 It is easy to check that
 $  (4x(x-1)(x-t))^{n/2}  f(\phi(x))$ satisfies
\begin{eqnarray*} p(x)y''+ (\frac{1}{2}-n) p'(x)y'+(4(H-n^2(t+1)+4n(2n-1)x)y&=&0, 
\end{eqnarray*}
which is the Halphen transform of the Lam\'e equation.
\end{proof}

\section{Examples}\label{Ex}

 In this section we go the opposite way of  Corollary~\ref{Mu}
to obtain
 a list of Lam\'e equations with arithmetic Fuchsian monodromy group.
 We start with a  list of 13 second order differential equations
 having $4$ regular singularities (Heun equations). 
 All these  arise from rational pull-backs of the Gauss hypergeometric function
 ${}_2F_1(\frac{1}{12},\frac{1}{12},\frac{2}{3},x).$ 
 In addition after scaling their monodromy group is contained in $\SL_2(\ZZ)$ 
 and they  posses at least $3$ unipotent monodromy group
 generators. 
 Some of the 17 differential equations 
 we obtained in \cite{MR08-2} from Herfurtner's list, where
 Herfurtner has classified families of elliptic curves with $4$ 
 singular fibres, 
 coincide. Thus only 13 remain.

\begin{tab}\label{Her1-18}
 Heun equations
 having at least $3$ unipotent
 monodromy group generators and monodromy group in $\SL_2(\ZZ)$
 being pull-backs of  Gauss hypergeometric differential equations,
 taken from \cite{MR08-2}: 
\[
\begin{array}{llll}
   \mbox{Nr. in} & \mbox{Herfurtner's} &  p(x)y''+p'(x)y'+(\alpha \beta x+\tilde{H})y=0& p(x)\\
  \cite{MR08-2}   &  \mbox{Notation} &&\\
  1 &I_1 I_1 I_1 I_9 & p(x)y''+p'(x)y'+(x+1/3)y=0 &x(x^2+x+1/3) \\
 2 &I_1 I_1 I_2 I_8 &  p(x)y''+p'(x)y'+xy=0 &  x(x-1)(x+1) \\
 3 &I_1 I_2 I_3 I_6&  p(x)y''+p'(x)y'+(x-1/4)y=0& x(x-1)(x+1/8)\\
  4& I_1 I_1 I_5 I_5& p(x)y''+p'(x)y'+(x+3)y =0 &x(x^2+11x-1) \\
  7& I_1 I_1 I_8 II   & p(x)y''+p'(x)y'+(35/36x - 3/2)y =0 &  x(x^2-14/3x+9) \\
 8 &I_1 I_2 I_7 II&  p(x)y''+p'(x)y'+(35/36x-9/8)y=0&   x(x-1)(x-81/32) \\
 9 &I_1 I_4 I_5 II &  p(x)y''+p'(x)y'+(35/36x-45/2)y=0,&   x(x-1)(x-81) \\
 10& I_2 I_3 I_5 II &  p(x)y''+p'(x)y'+(35/36x-20/27)y=0,&   x(x-1)(x-32/27) \\
  11& I_1 I_1 I_7 III  & p(x)y''+p'(x)y'+(15/16x +1/8)y=0 & x(x^2+13/32x+1/8)\\
12& I_1 I_2 I_6 III &  p(x)y''+p'(x)y'+(15/16 x-3/2)y=0,&   x(x-1)(x-4) \\
 13&  I_1 I_3 I_5 III;&  p(x)y''+p'(x)y'+(15/16 x-12)y=0,&   x(x-1)(x-128/3) \\
 15 & I_1 I_1 I_6 IV&  p(x)y''+p'(x)y'+(8/9 x)y=0,&   x(x-1)(x+1) \\
 16 &I_1 I_2 I_5 IV&  p(x)y''+p'(x)y'+(8/9x-8/27)y=0,&   x(x-1)(x-2/27) \\
\end{array}
\]

Herfurtner's list can be seen to be  indexed by the local monodromy in $\SL_2(\ZZ)$,
where $I_k$ denotes the  unipotent class in $\SL_2(\ZZ)$ containing a
 triangular matrix with entry $k$ in the upper diagonal and $II, III$ and $IV$ classes of 
 elliptic elements
 of order $6, 4$ and $3$ resp.
\end{tab}

This list gives rise to the following list of Lam\'e equations via
Corollary~\ref{Mu} a)
and the relations between the coefficients of the Heun equation and the Lam\'e
equation given there:
(Note that $p(x)$ changes to $4p(x)$.)

\[ -n(n+1)=\alpha \beta -3/4, \quad  H=\tilde{H}+\frac{1}{4} \sum_{i=1}^3 e_i\]
Thus we get the following table of geometric Lam\'e equations:

\begin{tab}
\[
\begin{array}{lllll}
   \mbox{Nr.} &
  \mbox{\cite[p.23]{Chud1383}} &  p(x)y''+\frac{1}{2}p'(x)y'- ( n(n+1)x-H)y=0& p(x)=\\
 1 & & p(x)y''+\frac{1}{2}p'(x)y'+(1/4x+1/12)y=0 &4x(x^2+x+1/3) \\
 2 & &  p(x)y''+\frac{1}{2}p'(x)y'+1/4xy=0 &  4x(x-1)(x+1) \\
 3 & &  p(x)y''+\frac{1}{2}p'(x)y'+(1/4 x-1/32)y=0& 4x(x-1)(x+1/8)\\
  4&  & p(x)y''+\frac{1}{2}p'(x)y'+(1/4x+1/4)y=0  &4x(x^2+11x-1) \\
 7 &    & p(x)y''+\frac{1}{2}p'(x)y'+(2/9x - 1/3)y =0 &  4x(x^2-14/3x+9) \\
 8 & (1,3);(2) &  p(x)y''+\frac{1}{2}p'(x)y'+(2/9x-31/128)y=0&  4 x(x-1)(x-81/32) \\
 9 &(1,3); (4) &  p(x)y''+\frac{1}{2}p'(x)y'+(2/9x-2)y=0,&  4 x(x-1)(x-81) \\
 10&    &  p(x)y''+\frac{1}{2}p'(x)y'+(2/9x-7/36)y=0,& 4  x(x-1)(x-32/27) \\
 11&   & p(x)y''+\frac{1}{2}p'(x)y'+(3/16x +3/128)y=0 & 4x(x^2+13/32x+1/8)\\
 12 & (1,2);(2) &  p(x)y''+\frac{1}{2}p'(x)y'+(3/16 x-1/4)y=0,&  4 x(x-1)(x-4) \\
 13  & (1,2);(3)&  p(x)y''+\frac{1}{2}p'(x)y'+(3/16 x-13/12)y=0,&  4 x(x-1)(x-128/3) \\
 15 &  &  p(x)y''+\frac{1}{2}p'(x)y'+(5/36 x)y=0,&  4 x(x-1)(x+1) \\
 16 & &  p(x)y''+\frac{1}{2}p'(x)y'+(5/36x-1/36)y=0,& 4  x(x-1)(x-2/27) \\
\end{array}
\]
The entry $(1,e);(\cdot)$ in the second column refers to the Lam\'e equation $(\cdot)$
  with arithmetic
 monodromy group of signature $(1,e)$ given in
\cite[p. 23]{Chud1383}.
\end{tab}

\begin{rem}
 \begin{enumerate}
 \item
 It is mentioned in \cite[Section 3]{Chud1383} that for all $71$
 (s. \cite{take83})
  arithmetic Fuchsian subgroups $\Gamma$ of signature $(1,e)$
 \[ \Gamma=\langle \alpha ,\beta, \gamma \mid   \alpha \beta \alpha^{-1} \beta^{-1} \gamma=-1_2, \gamma^e=-1_2\rangle ,\]
 where $\alpha$ and $\beta$ are hyperbolic elements of $\SL_2(\RR)$,
 there exists a corresponding Lam\'e equation,
 defined over $\overline{\QQ}$.
  Using numerical solutions of the (inverse)
 uniformization problem for the punctured tori
 they were computed. But only some
 of these were listed there.
 Our list of Lam\'e equations contains all those where the 
 quaternion algebra $A$ over $k$ associated to the arithmetic Fuchsian group
 is a  quaternion algebra $A$ over $\QQ$ (s. \cite{take83}), i.e.
 Nr. 1,2,3,4, 7,8,9,10,11,12,13.
 \item
 The examples 15 and 16 do not have an arithmetic monodromy group of signature $(1,e)$.
  However since
 \[ |\tr (A_1A_2A_3)|=1=2\cos(\pi/3) ,\]
the monodromy group is an arithmetic Fuchsian group with signature $(0;2,2,2,3),$
 s. \cite[Thm. 1]{Pu}.
 \end{enumerate} 
\end{rem}

\begin{proof}
 It follows from the following section that the monodromy group of the
 Lam\'e equation (in $\SL$-form) is contained
 in $\SL_2(\RR)$.
 Thus the comparison with the Fricke parameters in
 \cite{take83} (in \cite{ANR00} resp.)
 and Lemma~\ref{Min-Tup} (using Lemma~\ref{trace}) yields the claim. 
\end{proof}

Lam\'e equations with unipotent monodromy at infinity were already studied by Chudnovsky and Chudnovsky
in \cite{Chud1383} via the Halphen transform and symmetric squares.
There it was also mentioned that the Heun cases 1,2,3,4 are pull-backs of hypergeometric
differential equations. Using computer aided computations the following conjecture was stated: 

\begin{conj}[Chudnovsky and Chudnovsky]
 Lam\'e equations with $n=-\frac{1}{2}$ defined over $\bar{\QQ}$ are not
 globally nilpotent except for the 4 classes listed as $1,2,3,4$ above. 
\end{conj}

Note that Beukers also
studied  Heun equations with $4$ unipotent monodromy group generators
in \cite{Beukers02} and Lam\'e equations with unipotent monodromy
in \cite{Beukers07}.

Next we list Heun equations with $3$ unipotent
monodromy group generators obtained via rational Belyi functions (i.e. rational functions which are only ramified at $0,1$ and $\infty$)  
that do not appear
in \cite{MR08-2}. In these cases the monodromy group is a subgroup
of a nonarithmetic triangle group.
The conditions for right choice of Belyi-functions $j(x)$
and the hypergeometric differential equation  follow from \cite[Sec. Belyi functions]{MR08-2}. Since the computation of the Heun equations
is analogous to the one in \cite{MR08-2} we skip it.  
We only list  $j(x)$ and the corresponding hypergeometric function that
yield the Heun equation.

\begin{lemma}
Let $j(x)=\frac{j_1(x)}{j_2(x)}$ and   ${}_2 F_1 (a,b,c,x)$ be as in the list below:
\[
    \begin{array}{lccll}
     &    j(x)  &  \mbox {ramification data } &{}_2 F_1 (a,b,c,x)  & \mbox{ Riemann scheme} \\
   i)&   -\frac{(x^2-10x+5)^2x}{(5x^2-10x+1)^2} & (2,2,1),(5),(2,2,1)& 
{}_2 F_1 (\frac{13}{20},\frac{3}{20},1,x) & \left(\begin{array}{ccc}  
                           0 & 1 & \infty \\
                            0 & 0 & \frac{13}{20} \\
                            0 & \frac{1}{5} & \frac{3}{20}
                           \end{array}\right).\\
ii)& \frac{64(3x-1)^5}{27 x^3 ( 576x^2-285x+40)}&(5),(2,2,1),(3,1,1)&
 {}_2 F_1(\frac{3}{20},\frac{3}{20},\frac{4}{5},x)&\left(\begin{array}{ccc}  
                           0 & 1 & \infty \\
                            0 & 0 &  \frac{3}{20} \\
                            \frac{1}{5} & \frac{1}{2} & \frac{3}{20}
                           \end{array}\right).\\
iii)& \frac{(x+80)^3x^2(25x-48)}{64 (3x-16)^5}&(3,2,1),(2,2,2),(5,1)&
   {}_2 F_1(\frac{7}{20},\frac{3}{20},1,x)&\left(\begin{array}{ccc}  
                           0 & 1 & \infty \\
                            0 & 0 &  \frac{7}{20} \\
                            0 & \frac{1}{2} & \frac{3}{20}
                           \end{array}\right).\\   
 iv)& -\frac{x^4(25x^2+44x+20)}{256(x+1)^5}& (4,1,1),(2,2,2),(5,1)&   {}_2 F_1(\frac{7}{20},\frac{3}{20},1,x)
 &\left(\begin{array}{ccc}  
                           0 & 1 & \infty \\
                            0 & 0 &  \frac{7}{20} \\
                            0 & \frac{1}{2} & \frac{3}{20}
                           \end{array}\right).\\

 v)&-\frac{x^2(81x^2+176x+96)}{256(x-1)^3}& (2,1,1),(4),(3,1)&
 {}_2 F_1(\frac{13}{24},\frac{5}{24},1,x)&\left(\begin{array}{ccc}  
                           0 & 1 & \infty \\
                            0 & 0 &  \frac{13}{24} \\
                            0 & \frac{1}{4} & \frac{5}{24}
                           \end{array}\right).\\
\end{array}\]

Then the function $j_2(x)^{-a}\; {}_2 F_1(a,b,c,j(x))$ is a solution of

\[  \begin{array}{lllcc}
     &  p(x)y''+p'(x)y'+(\alpha \beta x-q)y=0 & p(x) \\
 i) & p(x)y''+p'(x)y'+(15/16x -9/16)y=0 & x(x^2-2x+1/5) \\
 ii) &p(x)y''+p'(x)y'+(15/16x-9/4)y=0& x(x^2-57/8x+72/5) \\
 iii) & p(x)y''+p'(x)y'+(99/100x+45/4)y=0    &x(x-1)(x+125/3)\\
 iv)  & p(x)y''+p'(x)y'+(99/100x+3/4)y=0& x(x^2+11/5x+5/4) \\
 v)  &p(x)y''+p'(x)y'+( 35/36 x+20/27)y=0& x(x^2+176/81x+32/27)\\
   \end{array}\]
\end{lemma}

\begin{cor}
The examples in the above Lemma provide the following Lam\'e equations:
\[  \begin{array}{cllcc}
     &  p(x)y''+\frac{1}{2}p'(x)y'-(n(n+1)x-H))y=0 & p(x) \\
 i) & p(x)y''+\frac{1}{2}p'(x)y'+(3/16x -1/16)y=0 & 4x(x^2-2x+1/5) \\
 ii) &p(x)y''+\frac{1}{2}p'(x)y'+(3/16x-15/32)y=0& 4x(x^2-57/8x+72/5) \\
 iii) & p(x)y''+\frac{1}{2}p'(x)y'+(6/25x+13/12)y=0    &4x(x-1)(x+125/3)\\
 iv)  & p(x)y''+\frac{1}{2}p'(x)y'+(6/25x+1/5)y=0& 4x(x^2+11/5x+5/4) \\
 v)  &p(x)y''+\frac{1}{2}p'(x)y'+( 2/9 x+16/81)y=0& 4x(x^2+176/81x+96/81)\\
   \end{array}\]
Example iii) appears (after the M\"obius transformation
$\phi(x)= (-125/3-1)x+1$) in \cite[p.23/24]{Chud1383} and has an arithmetic
Fuchsian monodromy group with signature $(1,5)$.
\end{cor}

\begin{proof}
This follows from Corollary~\ref{Mu} as in the computation of Table 3.2.
\end{proof}

Also geometric Heun equations with $4$ equal exponent differences and monodromy group contained in $\SL_2(\RR)$ yield after taking
the inverse Halphen transform geometric Lam\'e equations:
Those Heun equations can be computed for example as rational pull-backs of
Gauss hypergeometric differential equations with local projective
monodromy of orders
$(2,3,7), (2,3,8),$
$(2,3,9), (2,3,10)$ or $(3,3,5).$ 
In this cases the Heun equation has
local projective monodromy orders $(3,3,3,3).$
We demonstrate this via the following example:

\begin{ex}
 If one uses the Belyi function $j(x)=-\frac{(x-1)^3(x^2+3x+6)}{10x^2-15x+6}$
 with ramification data $(3)(1)(1),\; (5),\;$ $(3)(1)(1)$ as pull-back for
 ${}_2 F_1 (\frac{1}{15},\frac{2}{5},\frac{2}{3},x)$ one obtains after
 a M\"obius transformation the Heun equation
 \[p(x)y''+\frac{2}{3} p'(x) y'+(2/9x-44/243)y=0,\; p(x)=x(x-1)(x-32/27). \]
 The corresponding monodromy group is a subgroup of finite index of
 the arithmetic triangle group corresponding to ${}_2 F_1 (\frac{1}{15},\frac{2}{5},\frac{2}{3},x)$.
 Applying the inverse Halphen transformation we get using Corollary~\ref{Mu} c)
 the Lam\'e equation
 \[ p(x)y''+\frac{1}{2} p'(x) y'-( n(n+1)x-H)y=0,\quad p(x)=4x(x-1)(x-32/27),\quad
   n=-\frac{1}{6}, H=-\frac{13}{108}.\] 

\end{ex}

\section{Monodromy}\label{Mono}

Here we determine the monodromy group generators of the Heun equations in
Table 3.1.
We will state some remarks concerning the change of the monodromy group generators
under the symmetric square and the Euler-integral with $\mu=-\frac{1}{2}$.
This allows us to determine the monodromy group generators of the corresponding
Lam\'e equations.

\begin{defn}
 We call a tuple $\A=(A_1,\ldots,A_4), A_i\in \SL_2(\CC),$ 
 a tuple of monodromy group generators in $\SL_2(\CC)^4$ if
 \[ A_1 A_2 A_3 A_4 =\id_2 \]
 and the eigenvalues of the $A_1,\ldots,A_4$ are roots of unity.
\end{defn}

It is well known that the monodromy group representation is uniquely
determined by the Fricke parameters:

\begin{thm}\cite[p. 365-366]{FrKl}
 Let  $\A$ be a tuple of monodromy group generators in $\SL_2(\CC)^4$,
\[a_1=\tr(A_1),\quad a_2=\tr(A_2),\quad  a_3=\tr(A_3),\quad  a_4=\tr(A_4)\]
 and 
 \[ x=\tr(A_1A_2),\quad y=\tr(A_2A_3),\quad z=\tr(A_1A_3). \]
 Then the parameters $(a_1,a_2,a_3,a_4,x,y,z)$ satisfy the Fricke relation
\[ \sum_{i=1}^4 a_i^2 +\prod_{i=1}^4a_i+
x^2+y^2+z^2+xyz-(a_1a_2+a_3a_4)x- (a_1a_4+a_2a_3)y-(a_1a_3+a_2a_4)z=4.\]
\end{thm}

A nice well known application is the following

\begin{cor}\label{form}
 Let the monodromy group act irreducibly.
  Then
 it leaves a hermitian form invariant if and only if
 all Fricke parameters are real numbers.

 If the form is positive definite then the group is contained
 in $\SU_2(\RR)$ and if it is indefinite then the group is contained
 in $\SL_2(\RR)$.
\end{cor}

\begin{cor}\label{Lame-form}
 Let the monodromy group of a Lam\'e equation  act irreducibly
 and leave an indefinite hermitian form invariant. 
 Then the  Fricke parameters $(x,y,z)$
 are  of absolute value $\geq 2$. 
\end{cor}

\begin{proof}
 The
 group generated by $A_1$ and $A_2$ is an irreducible dihedral group
 for $x\neq \pm 2$.
 If the product would be an elliptic element, i.e. $\mid x \mid <2$, 
 then
  this subgroup would leave a positive definite form invariant.
 Thus the claim follows.  
\end{proof}

Since the monodromy group is invariant under the action of the braid group
we at first consider the braid group orbit.
It is quite natural to take as representative a monodromy tuple
 with minimal Fricke parameters. Thus we recall

\begin{lemma}
 The braid group $B_2=\langle \beta_1, \beta_2 \mid \beta_1 \beta_2 \beta_1=\beta_2 \beta_1 \beta_2\rangle$ acts
 on $\A\in \SL_2(\CC)^4$ via
 \begin{eqnarray*} \beta_1(\A)=(A_2,A_2^{-1} A_1 A_2 ,A_3,A_4)&&
     \beta_2(\A)=(A_1,A_3,A_3^{-1} A_2 A_3, A_4).   \end{eqnarray*}
 This yields the following transformation of the Fricke-parameters:
\[ \beta_1: (a_1,a_2,a_3,a_4,x,y,z)\mapsto (a_2,a_1,a_3,a_4,x,\tilde{z},y),\quad
  \tilde{z}=a_1a_3+a_2a_4-z-xy,\]
 \[ \beta_2: (a_1,a_2,a_3,a_4,x,y,z)\mapsto (a_1,a_3,a_2,a_4,z,y,\tilde {x}),\quad
  \tilde{x}=a_1a_2+a_3a_4-x-yz.\]
\end{lemma}

We consider the case where we have  local unipotent monodromy
at least at $3$ singularities. 

\begin{cor}\label{Fricke}
  Let $A_1,\ldots, A_3$ be unipotent elements
and
\[ x=n_1+2,\quad  y=n_2+2,\quad z=n_3+2.\] 
Then the Fricke relation reads
 \begin{eqnarray*}(\sum n_i+2-a_4)^2+\prod  n_i&=&0. \end{eqnarray*}
 Further the second solution $n_i'$
 of the quadratic equation for $n_i$ 
 is obtained via the corresponding braid group action  
 and we get $n_in_i'=(n_j+n_k+2-a_4)^2.$
Moreover, for the special value of $a_4=2$
 the Fricke relation simplifies to
\begin{eqnarray*} (x-4)^2+(y-4)^2+(z-4)^2&=&20-xyz. \end{eqnarray*}
If we put in the case $a_4=2$
\[ x=n_1N+2,\quad   y=n_2N+2,\quad  z=n_3N+2,\quad N=gcd(n_1,n_2,n_3) \in \NN\]
 we obtain
 \begin{eqnarray*} (n_1+n_2+n_3)^2+n_1n_2n_3N &=&0. \end{eqnarray*}
\end{cor}

\begin{proof}
 This follows from the above lemma 
 and direct computations with the Fricke relation
 using the identities:
 \[ x^2+y^2+z^2+xyz-4=(\sum n_i+4)^2+\prod n_i, \]
\[ \sum_{i=1}^4 a_i^2 +\prod_{i=1}^4a_i
     -(a_1a_2+a_3a_4)x- (a_1a_4+a_2a_3)y-(a_1a_3+a_2a_4)z=
   (2+a_4)^2-2(2+a_4)(\sum n_i+4).\]
 Since  $a_1=a_2=a_3=2$ the claim is readily to check.
\end{proof}

\begin{cor}\label{notgeom}
 Let $\A\in \SL_2(\CC)^4$ be a tuple of monodromy group generators 
 with unipotent
 elements $A_1,\ldots,A_3$, where one of the Fricke parameters is greater
 than $2$.
  Then in the braid group orbit
  there exists a tuple with Fricke parameters
  \begin{eqnarray*} (x,y,z) &=\left\{ \begin{array}{ccc}
                                    (x,2,-x+2+a_4)& x\geq 2 \\
                                  (x,-2,x)& x\geq 2,& a_4=-2.
                                 \end{array}\right.\end{eqnarray*}

The corresponding  monodromy group generators are
 \[\begin{array}{ccccc}

 \left(\begin{array}{cc}
       1 & 1 \\
       0 & 1 
        \end{array} \right),&
 \left(\begin{array}{cc}
       1 & 0 \\
       x-2 & 1 
        \end{array}\right) ,&
 \left(\begin{array}{cc}
       1 & 0 \\
       -x+a_4 & 1 
        \end{array}\right) ,&
 \left(\begin{array}{cc}
       1 & -1 \\
       -a_2+2 & a_4-1 
        \end{array}\right) &\\
 \left(\begin{array}{cc}
       1 & x-2 \\
       0 & 1 
        \end{array}\right) ,&
 \left(\begin{array}{cc}
       3 & -4 \\
       1 & -1 
        \end{array}\right) ,&
 \left(\begin{array}{cc}
       1 & 0  \\
        1 & 1
        \end{array}\right) ,&
 \left(\begin{array}{cc}
       -1 & x+2  \\
       0 & -1 
        \end{array}\right).
\end{array}\]
\end{cor}

\begin{proof}
 Let $x=n_1+2, y=n_2+2, z=n_3+2$.
 By Corollary~\ref{Fricke} we have
 \[ (\sum n_i+a)^2+\prod  n_i=0,\quad  0\leq 2-a_4=a\leq 4. \]
 Via braiding we can assume therefore that there are Fricke parameters $(x,y,z)$
 satisfying
 $x\geq y \geq 2\geq z$.
 We choose such a triple $(x,y,z)$, where $y\geq 2 $ is minimal.
 Further we take the minimal
 $x$ under the action of $\langle \beta_2^2 \rangle.$
 Then  by Corollary~\ref{Fricke} the minimality implies
 \[ |n_2+n_3+a|\geq 
 n_1\geq n_2\geq 0 \geq n_3 \geq -(n_1+n_2+a).  \]

 Case i) Let  $n_1+n_3+a \geq 0.$
 If $n_3+a\geq 0$ then
 \[ (n_1+n_2)^2 \leq (n_1+n_2+n_3+a)^2 =-n_1n_2n_3 \leq n_1n_2a \leq 4 n_1n_2.\] This implies
 \[ n_1=n_2,\quad a=4=-n_3. \]
 If $n_3+a< 0$ then
  \[ n_1 \leq |n_2+n_3+a| \leq |n_3+a| \leq n_1. \]
 Hence $n_2=0$ and therefore $n_1+n_3+a=0$.

 Case ii) Let  $n_1+n_3+a < 0.$
 Then
  \[ n_2^2 \geq ((n_1+n_3+a)+n_2)^2 =-n_1n_2n_3>n_1n_2(n_1+a)>n_1^2n_2, \]
 since  $-n_3\leq (n_1+n_2+a).$
 Thus $n_2=0$.
\end{proof}

\begin{lemma}\label{negative}
 Let $x, y, z$ be negative integers satisfy the
 relation in Corollary~\ref{Fricke}. Then
 there exist in the braid group orbit the following minimal 
 triple (with respect to $|x|\leq |y|\leq |z|$):

\[ \begin{array}{lllllll}
\mbox{ case}&  a_4 & & (n_1,n_2,n_3)&\\ 
i)& a_4=2&N=5 &(-1,-4,-5)\\
&&N=6 &(-1,-2,-3)\\
&&N=8 &(-1,-1,-2)\\
&&N=9 &(-1,-1,-1)\\                                   
ii)&   a_4=0 & & (-5,-12,-15) &(-6-8,-12)&(-7,-7,-9)\\
iii)&   a_4=1&& (-5,-16,-20)&(-6,-10,-15) &(-7,-8,-14)&\\
    &&&(-8,-8,-9)\\
iv)&   a_4=-1 &&(-5,-8,-10)&(-6,-6,-9)
   \end{array} \]
\end{lemma}

\begin{proof}
Case i):  
Let $a_2=2$ and
  $n_3\leq n_2 \leq n_1<0$ be a minimal triple. Thus  $-n_3\leq -(n_1+n_2)$
and therefore
\begin{eqnarray*} 4(n_1+n_2)^2\geq (n_1+n_2+n_3)^2=-Nn_1n_2n_3.\end{eqnarray*}
Hence
\begin{eqnarray*} 4(\frac{-1}{n_1}+\frac{-1}{n_2})^2 \geq \frac{-Nn_3}{n_1n_2}
= \frac{(Nn_3)^2}{-n_1n_2n_3N}\geq \frac{(Nn_3)^2}{(3n_3)^2} \geq \frac{N^2}{9}\end{eqnarray*}
and we obtain
\begin{eqnarray*} (\frac{-1}{n_1}+\frac{-1}{n_2})\geq \frac{N}{6}.\end{eqnarray*}

At first we consider solutions with $n_1=-1$:
Then   $n_3\in \{n_2,n_2-1\}$.
If the tuple is of the form
\begin{enumerate}
\item[i)] $(-1,n_2,n_2)$ then
\begin{eqnarray*} (-1-n_2-n_2)^2=Nn_2^2 \Rightarrow n_2 =-1. \end{eqnarray*}
Thus $N=9$ and $(-1,n_2,n_2)=(-1,-1,-1)$.

\item[ii)] $(-1,n_2,n_2-1)$ then
\begin{eqnarray*} (-1+n_2+n_2-1)^2=Nn_2(n_2-1) \Rightarrow 4(n_2-1) =Nn_2. \end{eqnarray*}
Hence $N>4$ and $n_2 \mid 4$.
This gives the cases
 $N=8$ and $(-1,n_2,n_2-1)=(-1,-1,-2)$
or $N=6$ and $(-1,n_2,n_2-1)=(-1,-2,-3)$
or $N=5$ and $(-1,n_2,n_2-1)=(-1,-4,-5)$.
\end{enumerate}

All other solutions start with at least $n_1=-2$. Hence
$N\leq 6.$

Next we consider solutions with $n_1=-2$:
Then  $n_3\in \{n_2,n_2-1,n_2-2\}$.
If the tuple is of the form
\begin{enumerate}
\item[i)] $(-2,n_2,n_2)$ then
\begin{eqnarray*} (-2-n_2-n_2)^2=2Nn_2^2 \Rightarrow n_2 \mid 2. \end{eqnarray*}But $n_2=2$ yields a contradiction.

\item[ii)]  $(-2,n_2,n_2-1)$ we get
\begin{eqnarray*} (-2+n_2+n_2-1)^2=Nn_2(n_2-1) \Rightarrow n_2-1 \mid 1,\quad n_2\mid 3 \end{eqnarray*}
a contradiction.

\item[iii)] $(-2,n_2,n_2-2)$ then
\begin{eqnarray*} (-2+n_2+n_2-2)^2=Nn_2(n_2-2) \Rightarrow 2(n_2-2)=Nn_2.\end{eqnarray*}
Thus $n_2\mid 4$ and $N-2\mid 4.$
This  yields the cases
\begin{eqnarray*} N=4,& (-2,-2,-4),\\
                   N=3, & (-2,-4,-6).\end{eqnarray*}
\end{enumerate}

For the remaining cases we only have to check  $N\leq 4.$

For $N=4$ we get the solution $(-3,-3,-6)$.
For $N=3$ we get $(-3,-3,-3)$ and
for  $N=2$ we get $(-3,-6,-9)$ or
$(-4,-4,-8)$.
But in these cases
$gcd(n_1,n_2,n_3)>1.$

Finally assume that $1=N$.
Then
\begin{eqnarray*} (n_1+n_2+n_3)^2=-n_1n_2n_3 \end{eqnarray*}
implies that $gcd(n_i,n_j)=1$ for $i \neq j.$
Hence $-n_i$ is a square.
But the equation  reduced modulo $3$ has no solutions.

For the other cases the proof is analogous.
\end{proof}

 For case i) see also Gutzwiller \cite{Gutz87}.

Now we can easily determine the corresponding tuple
of monodromy group generators:

\begin{lemma}\label{Min-Tup}
We list the corresponding tuples of monodromy group generators for
the minimal tuples in Lemma~\ref{negative}.
The monodromy group generators for  $N=3$ and $N=9$
are conjugate in $\SL_2(\QQ)$ by a diagonal matrix. The same holds for 
 $N=4$ and $N=8$.

\[
\begin{array}{ccccccc}
&&A_1& \quad A_2& \quad A_3& \quad A_4\\
&N=3,&  \left(\begin{array}{cc}
      1 & 0 \\
      -3 & 1 
         \end{array} \right),& 
  \left(\begin{array}{cc}
      -2 & 3\\
      -3 & 4 
         \end{array} \right),& 
 \left(\begin{array}{cc}
      -5 & 12\\
      -3 & 7 
         \end{array} \right),& 
 \left(\begin{array}{cc}
      1 & 3\\
      0 & 1 
         \end{array} \right) \\
&N=4,&  \left(\begin{array}{cc}
      1 & -2\\
      0 & 1 
         \end{array} \right),& 
  \left(\begin{array}{cc}
      5 & -2\\
      8 & -3 
         \end{array} \right),& 
 \left(\begin{array}{cc}
      1 & 0\\
      4 & 1 
         \end{array} \right),& 
 \left(\begin{array}{cc}
      -3 & -4\\
      4 & 5 
         \end{array} \right) \\
&N=5,&
 \left(\begin{array}{cc}
      1 & 5\\
      0 & 1 
         \end{array} \right),& 
  \left(\begin{array}{cc}
      1 & 0\\
      -1 & 1 
         \end{array} \right),& 
 \left(\begin{array}{cc}
      -9 & 20\\
      -5 & 11 
         \end{array} \right),& 
 \left(\begin{array}{cc}
      -9 & 25\\
      -4 & 11 
         \end{array} \right)\\
&N=6,&  \left(\begin{array}{cc}
      1 & 0 \\
      -6 & 1 
         \end{array} \right),& 
  \left(\begin{array}{cc}
      -5 & 2\\
      -18 & 7 
         \end{array} \right),& 
 \left(\begin{array}{cc}
      -5 & 3\\
      -12 & 7 
         \end{array} \right),& 
 \left(\begin{array}{cc}
      1 & 1\\
      0 & 1 
         \end{array} \right) \\
&N=8,&  \left(\begin{array}{cc}
      1 & -1\\
      0 & 1 
         \end{array} \right),& 
  \left(\begin{array}{cc}
      5 & -1\\
      16 & -3 
         \end{array} \right),& 
 \left(\begin{array}{cc}
      1 & 0\\
      8 & 1 
         \end{array} \right),& 
 \left(\begin{array}{cc}
      -3 & -2\\
      8 & 5 
         \end{array} \right)\\
&N=9,&  \left(\begin{array}{cc}
      1 & 0 \\
      -9 & 1 
         \end{array} \right),& 
  \left(\begin{array}{cc}
      -2 & 1\\
      -9 & 4 
         \end{array} \right),& 
 \left(\begin{array}{cc}
      -5 & 4\\
      -9 & 7 
         \end{array} \right),& 
 \left(\begin{array}{cc}
      1 & 1\\
      0 & 1 
         \end{array} \right)\\
\end{array}
\]

\[
\begin{array}{ccccccc}
(n_1,n_2,n_3)&&A_1& \quad A_2& \quad A_3& \quad A_4\\
(-6,-6,-9)&&\left(\begin{array}{cc}
           1  &   6  \\
          0   &1
        \end{array}\right),&  
\left(\begin{array}{cc}
          1   &   0  \\
           -1  &1
        \end{array}\right),&  
\left(\begin{array}{cc}
          -2   &  9   \\
          -1   &4
        \end{array}\right),&  
\left(\begin{array}{cc}
          -5   &  21   \\
           -1  &4
        \end{array}\right)\\
(-5,-8,-10)&&\left(\begin{array}{cc}
         1    &   5  \\
         0    &1
        \end{array}\right),&  
\left(\begin{array}{cc}
          3   &   4  \\
          -1   & -1
        \end{array}\right),&  
\left(\begin{array}{cc}
          1   &  0   \\
         -2    &1
        \end{array}\right),&  
\left(\begin{array}{cc}
          -1   & 1    \\
          -1   & 0
        \end{array}\right)\\
(-5,-16,-20)&&\left(\begin{array}{cc}
         1    &   5  \\
         0    &1
        \end{array}\right),&  
\left(\begin{array}{cc}
          2   &   1  \\
          -1   &0
        \end{array}\right),&  
\left(\begin{array}{cc}
          -3   &  4   \\
          -4   &5
        \end{array}\right),&  
\left(\begin{array}{cc}
          -4   &7     \\
          -3   &5
        \end{array}\right)\\
(-6,-10,-15)& &\left(\begin{array}{cc}
          1   &   5  \\
          0   &1
        \end{array}\right),&  
\left(\begin{array}{cc}
          1   &  0   \\
           -2  &1
        \end{array}\right),&  
\left(\begin{array}{cc}
           -2  &   3  \\
           -3  &4
        \end{array}\right),&  
\left(\begin{array}{cc}
          -2   &   7  \\
           -1  &3
        \end{array}\right)\\
(-7,-8,-14)&&\left(\begin{array}{cc}
          1   &  7   \\
          0   &1
        \end{array}\right),&  
\left(\begin{array}{cc}
         1    &  0   \\
        -1     &1
        \end{array}\right),&  
\left(\begin{array}{cc}
          -3   &   8  \\
          -2   &5
        \end{array}\right),&  
\left(\begin{array}{cc}
         -3    &   13  \\
          -1   &4
        \end{array}\right)\\
(-8,-8,-9)&&\left(\begin{array}{cc}
          1   &  8   \\
          0   &1
        \end{array}\right),&  
\left(\begin{array}{cc}
          1   &   0  \\
           -1  &1
        \end{array}\right),&  
\left(\begin{array}{cc}
           -2  & 9    \\
           -1  &4
        \end{array}\right),&  
\left(\begin{array}{cc}
          -5   &31     \\
           -1  & 6
        \end{array}\right)\\
(-5,-12,-15)& &\left(\begin{array}{cc}
          1   &  5   \\
          0   &1
        \end{array}\right),&  
\left(\begin{array}{cc}
          1   &   0  \\
          -1   &1
        \end{array}\right),&  
\left(\begin{array}{cc}
          -5   &   12  \\
           -3  &7
        \end{array}\right),&  
\left(\begin{array}{cc}
           -5  &   13  \\
           -2  & 5
        \end{array}\right)\\
(-6-8,-12)&&\left(\begin{array}{cc}
           1  &  6   \\
          0   &1
        \end{array}\right),&  
\left(\begin{array}{cc}
          1   &   0  \\
          -1  &1
        \end{array}\right),&  
\left(\begin{array}{cc}
          -3   &  8   \\
          -2   & 5
        \end{array}\right),&  
\left(\begin{array}{cc}
          -3   &  10   \\
           -1  &   3
        \end{array}\right)\\
(-7,-7,-9)&&\left(\begin{array}{cc}
           1  &   7  \\
          0   &1
        \end{array}\right),&  
\left(\begin{array}{cc}
          1   &  0   \\
          -1   &1
        \end{array}\right),&  
\left(\begin{array}{cc}
           -2  & 9    \\
           -1  & 4
        \end{array}\right),&  
\left(\begin{array}{cc}
           -5  & 26     \\
           -1  &  5
        \end{array}\right)\\
\end{array}
\]

\end{lemma}

\begin{rem}
 These tuple correspond exactly (up to the braid group action)
 to the differential equations in Table~3.1.
\end{rem}

 Next we show how the inverse Halphen transform changes the Fricke parameters.
 Hence one can also determine explicitly from this list 
 the monodromy group generators
 of the corresponding Lam\'e equations.
 Moreover we can identify the arithmetic Fuchsian monodromy group of the Lam\'e equations in Table~3.2
 with those from Takeuchi's list in \cite[Thm. 4.1 ]{take83}.

\begin{thm}\label{trace}
 Let $\A$ be an irreducible tuple of monodromy group generators in $\GL_2(\CC)^4$
 with $A_1, A_2, A_3$ being involutions. 
 Taking the symmetric square $\sym^2$ and applying the middle convolution $MC_{-1}$
 (s. \cite{DR99}) we obtain 
 a tuple of monodromy group generators $\phi(\A):=MC_{-1}\circ \sym^2 (\A)=: \B$ in $\SL_2(\CC)^4$
 with unipotent elements $B_1, B_2, B_3.$
  
 The induced transformation of the Fricke parameters
 \begin{eqnarray*}
 \phi :\; \{(0,0,0, a_4,x,y,z)\} &\to&  \{(2,2,2,b_4,{x'},{y'},{z'})\},\\
   (0,0,0, a_4,x,y,z) &\mapsto& (2,2,2,-a_4^2-2,-(x^2-2),-(y^2-2),-(z^2-2))
  \end{eqnarray*} 
 is bijective (identifying $\A$ with $(\eps_1A_1,\ldots, \eps_4A_4),\;\eps_i \in \{\pm 1 \},\;\prod \eps_i=1$).
\end{thm}

\begin{proof}
 The symmetric square of $\A$ gives a tuple
 $\C$ in $\SO_3(\CC)^4$, where $C_1, C_2$ and $C_3$ are reflections.
 If $A_4$ is semi-simple with eigenvalues $\alpha,-\alpha^{-1}$ 
 then $C_4$ has eigenvalues $\alpha^2,\alpha^{-2},-1.$
 If $\alpha \neq \pm 1$ then $\C$ is irreducible. 
 Applying the middle convolution functor $MC_{-1}$ to $\C$ we get
 by \cite[Thm 2.4 i)]{DR07} a tuple
 \[  MC_{-1}(\C)=\B \in \GL_m(\CC)^4,\quad m=\sum_{i=1}^3 \rk(C_i-1)+\rk(-C_4-1)-3=2.\]
 The change of the eigenvalues under the middle convolution follows from
 \cite[Lemma 2.6]{DR07} and
  gives a tuple of monodromy group generators $\B$ in $\SL_2(\CC)^4$
  with unipotent elements $B_1, B_2, B_3$ and
  $B_4$ being semi-simple with eigenvalues $-\alpha^2,-\alpha^{-2}$.
  The eigenvalues of $C_i C_j$ are $\gamma_{ij},\gamma_{ij}^{-1},1$
 where $\tr(C_1C_2)=x^2-1,$ $\tr(C_2C_3)=y^2-1 $ and
  $\tr(C_1C_3)=z^2-1. $
   Since the convolution commutes also with coalescing (s. \cite[Lemma 5.6]{DR99})
   the  eigenvalues of $B_iB_j$ are $-\gamma_{ij},-\gamma_{ij}^{-1}$.
   Hence 
   \[ \tr(B_iB_j)=-\tr(C_iC_j)+1 \]
    and the claim follows.
If $\alpha = \pm 1$ then $\A$ generates an orthogonal group in $\GO_2(\CC)$ and
 therefore $\C$ is reducible. 
 Thus we can assume $\C \in \SO_2(\CC)^4$, where
 $C_1,\ldots, C_4$  are reflections.
 By \cite[Thm 2.4 i)]{DR07}  we get
 \[  MC_{-1}(\C)=\B \in \GL_m(\CC)^4,\quad m=\sum_{i=1}^3 \rk(C_i-1)+\rk(-C_4-1)-2=2\]
  with $3$ unipotent elements $B_1,\ldots,B_3$ and a negative unipotent element
  $B_4$.

 If $A_4$ is not semi-simple then its Jordan form is
 $i J(2),$ where $J(2)$ denotes a Jordan block of length $2$. 
 Then $C_4$ has Jordan form $-J(3).$
 By \cite[Thm. 2.4. i)]{DR07}  we get
 \[  MC_{-1}(\C)=\B \in \GL_2(\CC)^4\]
 with $4$ unipotent elements $B_i$.
 The claim  follows in this case as above.

 The multiplicativity of $MC$ (s. \cite[Thm. 2.4. ii]{DR07}) 
  implies that 
 $MC_{-1} (\B)=MC_{1}(\C)=\C,$
 where the last equality follows from \cite[Prop.3.2]{DR99}.
 Since $\sym^2\; \SL_2(\CC)= \SO_3(\CC)$
 the induced map on the Fricke parameters is bijective.

 \end{proof}

\section{Properties of geometric monodromy groups}\label{Hodge}

In this section we make use of some properties of
monodromy groups of geometric differential equations.
Using operations that preserve geometric differential equations,
like rational pull-backs, tensor products and middle convolution
we obtain criteria to rule out certain monodromy group generators
arising from geometric differential equations.

At first we collect some properties of the monodromy group of a geometric
differential equation:

\begin{thm}\label{Deligne}
 Let
 \[ \rho: \pi_1(S,x_0) \to \GL_n(\CC) \]
 be a monodromy representation of a Picard-Fuchs differential equation:
 Then
 \begin{enumerate}
 \item  \[ \rho: \pi_1(S,x_0) \to \GL_n(\ZZ) \]
        Moreover  the monodromy group is either symplectic or orthogonal.
\item   Each absolutely irreducible composition factor leaves a hermitian
        form invariant.
\item  The local monodromy is quasi-unipotent.

\item  To a given Riemann scheme there exist only finitely many
       such monodromy representations $\rho$.

\end{enumerate}
\end{thm}

\begin{proof}
 The monodromy group representation of a Picard-Fuchs equation
  arises from a variation of
 $\QQ$-Hodge-structures and leaves a $\ZZ$-lattice invariant (cf. \cite[0. Introduction]{Del}).
 Hence we get
 \[ \rho :\pi_1(S,x_0) \to \GL_n(\ZZ). \]
 Since these $\QQ$-Hodge-structures are polarized the monodromy group leaves a
 symmetric or antisymmetric form invariant.
 
By Prop. 1.13 in \cite{Del} the absolutely irreducible composition factors
are underlying a polarizable complex Hodge-structure and therefore preserve
 a hermitian form.

It is well known that the local monodromy is quasi-unipotent
(cf. \cite[Introduction]{Katz70}).
The last statement follows from Thm. 0.5 in \cite{Del}.  
\end{proof}

This theorem has the following consequences:

\begin{cor}\label{algint}
 An absolutely irreducible composition factor of a  monodromy representation of a Picard-Fuchs differential equation can be defined over an integral ring $R$
that is contained in an abelian extension of the splitting field of the representation.
\end{cor}

\begin{proof}
 This follows from Theorem 23.18 in \cite{CR} and the Remark following it.
\end{proof}

\begin{cor}\label{sign} 
 Let 
 \[ \rho:  \pi_1(S,x_0) \to \SO_n(\CC) \]
 be a monodromy group representation of a geometric differential equation.
 Then
 the monodromy group leaves a hermitian form invariant with signature
 $(p,n-p).$
 \begin{enumerate}
  \item  If  $n=4$
 then the  signature is in
 $\{ (p,4-p) \mid \; p \in \{ 0,2,4\} \}$.
\item If  $n=6$ then the  signature is in
 $\{ (p,6-p) \mid\; p \in \{ 0,2,3,4,6\}\}$.
\end{enumerate}
\end{cor}

\begin{proof}
 By Thm.~\ref{Deligne} the  monodromy group representation of a geometric differential equation leaves a hermitian form invariant.
 \begin{enumerate}
 \item  Since
 \[ \SO_4(\CC)=\SL_2(\CC) \ten \SL_2(\CC) \]
 the monodromy group representation $\rho$ arises from
 a tensor product of two $2$-dimensional 
 monodromy group   representations of
 geometric differential equations.
 Since each of these preserves a hermitian form with signature 
 $(p,2-p),\; p \in \{ 0,1,2\},$ the tensor product of the hermitian forms
 has signature in $\{ (4,0),(2,2),(0,4) \}.$
\item
 Since
 \[ \SO_6(\CC)=\Lambda^2 \SL_4(\CC) \]
 the monodromy group representation $\rho$ arises from
 a antisymmetric tensor product of a  $4$-dimensional 
 monodromy group representation  of a
 geometric differential equation.
 Since this preserves a hermitian form with signature
 $(p,4-p),\; p \in \{ 0,1,2,3\},$ the antisymmetric 
 tensor product of the hermitian form
 has signature in $\{ (6,0),(4,2),(3,3), (2,4), (0,6) \}.$
\end{enumerate}
\end{proof}

\begin{rem}
 The group
 \[ \SL_2(\CC) \hookrightarrow \SO_4(\CC),\quad A \mapsto A \ten \bar{A^{t}} \]
 preserves a hermitian form, since the traces of all elements are real numbers.  Further, the signature is $(3,1)$.
\end{rem}

Making use of exceptional isomorphisms of Lie groups in small
dimension  we derive the following:

\begin{cor}\label{nongeom}
  Let $g$ be an element of a monodromy group of a geometric
 differential equation that is  contained in $\SO_n(\CC)$.
\begin{enumerate}
  \item Then $\tr(g)\in \RR.$ 
 \item If $3\leq n\leq 6$ and
        $\rk(g-\id_n)=2$ then $\tr(g)\geq n-4$. 
\end{enumerate}.
\end{cor}

\begin{proof}
\begin{enumerate}
 \item Since $g$ preserves a hermitian form and $g$ is contained in a selfdual
     group the set of eigenvalues (with multiplicity) is invariant under the
     complex conjugation. Hence  we get $\tr(g)\in \RR.$
 \item If $g \in \SO_3(\CC)$ we find
 $g_1 \in \SL_2(\CC)$ such that
 the symmetric square of $g_1$ is $g$.
 Hence the eigenvalues of $g$ are squares. 
 Either the eigenvalues $\alpha,\alpha^{-1}$ of $g_1$
  have absolute value $1$ or they are real and $\mid \tr(g_1) \mid >2.$
 Both cases imply the claim.

 If $g \in \SO_4(\CC)$ we find
 $g_1,g_2 \in \SL_2(\CC)$ such that
 $g_1 \ten g_2=g$.
 The condition $\rk(g-\id_4)=2$ implies that the Jordan form of $g_1$ is same
 as the Jordan form of $g_2$.
 Hence the eigenvalues of $g$ are squares. 
 Either the eigenvalues $\alpha,\alpha^{-1}$ of $g_1$
  have absolute value $1$ or they are real and $\mid \tr(g_1) \mid >2.$
 Both cases imply the claim.
 
 If $g \in \SO_6(\CC)$ we find
 $g_1 \in \SL_4(\CC)$ such that
 $ \Lambda^2 g_1 =g$.
 Thus $\rk(g_1-\alpha^{\pm 1} \id_4) =2$ for some $\alpha \in \CC.$
 Since $\tr(g_1)\in \RR$ we get $\alpha \in \RR$ or $\bar{\alpha}=\alpha^{-1}.$
  Again both cases imply the claim.

 The case $g \in \SO_5(\CC)=  \Lambda^2 \Sp_4(\CC)$ is analogous to the above case.
\end{enumerate}

\end{proof}

\begin{cor}\label{Nongeom2}
Let $\B\in \SL_2(\CC)^4$ be a tuple of monodromy group generators that
arises from a geometric differential equation with Fricke parameters
$(2,2,2, a_4,x,y,z).$
Then $a_4,x,y,z$ are real algebraic numbers, $\mid a_4 \mid \leq 2$, such that
 ${x},{y},{z} \leq -2 $
or  $-2 \leq {x},{y},{z} \leq 2.$ 
Moreover all Galois conjugates
of $a_4,x,y,z$ are also real algebraic numbers that satisfy the conditions. 
\end{cor}

\begin{proof}
By Thm.~\ref{Deligne} we get $\mid a_4\mid \leq 2$ and
that $a_4,x,y,z$ are real algebraic numbers.
Let the Fricke parameters $x,y,z$ not satisfy the conditions
 ${x},{y},{z} \leq -2 $
or  $-2 \leq {x},{y},{z} \leq 2.$ 
Then using Thm.~\ref{trace} 
we obtain tuple $\A, \phi(\A)=\B,$ that has a non-real Fricke parameter.
Hence it preserves no hermitian form.
Thus again by  Thm.~\ref{Deligne} it is not geometric.
However $\A$ can be transformed to $\B$ via geometric operations.
This shows the first claim.
Since the Galois conjugate monodromy representation of a geometric
differential equation is again geometric by Thm.~\ref{Deligne} iii)
the second claim follows.

\end{proof}

\begin{rem}
 The braid groups orbits in Corollary~\ref{notgeom} are non geometric.
\end{rem}

\begin{rem}
 From Cor. 1.10 in \cite{Del} one knows  that the Fricke parameters
 $(x,y,z)$ are bounded for a geometric tuple of monodromy group generators 
 in $\SL_2(\CC)^4$. Therefore one could expect that in a  geometric braid group orbit
  that (at most) the minimal Fricke parameters  arise
 from a geometric differential equation.
\end{rem}

Next we sum up  some properties of the middle convolution
cf. \cite{katz96}, \cite{DR99}, \cite{DR07}
adapted to our situation.

\begin{lemma}\label{prop} 
 Let $\A$ be a tuple of monodromy group generators in $\SL_2(\CC)^4$.
 Then the middle convolution $MC_{-1}(\A)$ yields a tuple
 $(B_1,\ldots,B_4)$ in $\GO_{m}(\CC)^4$, where
 \[ m=\sum_{i=1}^3\rk(A_i-1)-\rk(-A_4-1). \]
 Further,
\begin{enumerate}
 \item
 if $A_i,i\leq 3,$ (resp. $-A_4$)  is unipotent then
 $B_i$ (resp. $B_4$)  is a reflection.
 If $(\alpha,\alpha^{-1}),\alpha \neq 1, $ are the eigenvalues of
 $A_i,i\leq 3,$ (resp. $-A_4$) then
  $(-\alpha,-\alpha^{-1}) $ are the non trivial eigenvalues of
 $B_i,i\leq 3,$ (resp. $B_4$).
 \item  If $(\alpha,\alpha^{-1})$ are the eigenvalues of $A_iA_j,(i,j) \in \{(1,2),(2,3),(1,3)\},$ then
   \[ \left\{\begin{array}{ccc}
           (-\alpha,-\alpha^{-1}) &  \rk(A_i-1)+\rk(A_j-1)=2\\
            (-\alpha,-\alpha^{-1},-1) & \rk(A_i-1)+\rk(A_j-1)=3\\
             (-\alpha,-\alpha^{-1},-1,-1) & \rk(A_i-1)+\rk(A_j-1)=4
       \end{array}\right. \]
are the non trivial eigenvalues of $B_iB_j$.

\end{enumerate}
 \end{lemma}

\begin{lemma}\label{MCFricke}
  Let 
  $\A \in \SL_2(\CC)^4$ with Fricke parameters
  $(\alpha_1+\alpha_1^{-1},\ldots,
  \alpha_4+\alpha_4^{-1},x,y,z)$.
  Let $\gamma^2=\prod \alpha_i^{-1}$
  and
  $\B=(\alpha_1 A_1,\ldots,\alpha_3 A_3,\gamma^2 \alpha_4 A_4)$
  then the tuple
  \[
    \C:= (\gamma \alpha_1MC_{\gamma^2}(\B)_1,\ldots, \gamma \alpha_3MC_{\gamma^2}(\B)_3, \gamma^{-1} \alpha_4MC_{\gamma^2}(\B)_4) \in \SL_2(\CC)^4 \]
  has the Fricke-parameters $(\alpha_1 \gamma+\alpha_1^{-1}\gamma^{-1},\ldots,
  \alpha_4\gamma+\alpha_4^{-1}\gamma^{-1},x,y,z)$.
\end{lemma}

As an application  we reprove the following theorem due
to Chudnovsky and Chudnovsky (cf. \cite[Thm. 2.1]{Chud1383}).

\begin{thm}\label{Chud}
 The monodromy group of an irreducible  Lam\'e equation that is generated by $4$ reflections
 is geometric if and only if it is a finite dihedral group.
 \end{thm}

\begin{proof}
 It is well known that the monodromy group of  an irreducible  Lam\'e equation  generated by $4$ reflections 
is a dihedral one.
If it arises from a Picard-Fuchs equation then all eigenvalues of an element in the
 monodromy group are algebraic by Corollary~\ref{algint}.
 
 If the preserved hermitian form is indefinite
 then by Corollary~\ref{Lame-form} one of the Fricke-parameters has absolute value greater than $2$.
 We can assume that $x<-2.$
 Applying the middle convolution $MC_{-1}$ we get a
 tuple of monodromy group generators $(A_1,A_2,A_3,A_4)$, 
  where $A_1,A_2,A_3$ are unipotent and $A_4$ is negative
 unipotent and $\tr(A_1A_2)=-x.$
 Multiplying $A_3$ and $A_4$ by $-1$ doesn't change $x$
 and by Lemma~\ref{prop} 
  $ MC_{-1} (A_1,A_2,-A_3,-A_4)$ yields a tuple
   $(B_1,\ldots,B_4)$ in $\GO_4(\CC)^4$,
   where $B_1$ and $B_2$ are reflections
  with the property that
  $\tr(B_1B_2)=\tr(B_3B_4)=x+2<0$ and
  $\rk(B_1B_2-1)=2=\rk(B_3B_4-1)$.
  A quadratic pull-back yields a tuple  $(B_3,B_4,B_3^{B_1},B_4^{B_1}) \in \SO_4(\CC)^4.$
  Thus by Corollary~\ref{nongeom} it
  arises from a tuple of monodromy group generators that leaves no
  hermitian form invariant. 

  Hence the monodromy group of the Lam\'e equation has to be definite.
  Thus we can assume that
  all Galois conjugates of the monodromy group are also definite.
  Therefore the monodromy group is an isometry group of a lattice.
  This implies that it is a finite dihedral group.
 \end{proof}

\section{Appendix}\label{Classification}

We classify all braid orbits in $\SL_2(\ZZ)^4$.
Since multiplying by $-1$ preserves $\SL_2(\ZZ)$ we only
classify up to this operation. We also omit the case
where one of the $A_i$ is a central element.

\begin{defn}
Let $\mid a_1 \mid,\ldots ,\mid a_4 \mid \leq 2$
and the Fricke parameters
$(a_1,a_2,a_3,a_4,x,y,z)$ be in the braid group orbit such that
$\mid x \mid $ is minimal.
Let further $y,z$ such that $\mid y \mid +\mid z\mid $ is minimal.
Then we call $(a_1,a_2,a_3,a_4,x,y,z)$ a minimal tuple (in the braid group
orbit).
We denote by $x_2=\max\{\mid y \mid ,\mid z \mid\}$ 
and  $x_1=\min\{\mid y \mid ,\mid z \mid\}$ and by $x_0=\mid x\mid.$
\end{defn}

In the following we derive some bounds for minimal tuples with integer
Fricke parameters.

The braid group action provides the following bound.

\begin{lemma}\label{min}
 Let  $(a_1,a_2,a_3,a_4,x,y,z)$ be a minimal tuple.
 Then we have the following inequality:
\[ x_2 \leq  \left\{ \begin{array}{cc}
         (x_0+x_1+4)  & x_0x_1\geq 6\\
          (x_0+x_1+5) & x_0=2, x_1=2  \\      
           (x_1+6) & x_0=1,x_1 \leq 6 \\
             (x_1+5) & x_0=0, 2 \leq x_1  \\
            7 & x_0=0, x_1 \leq 1  \\
           \end{array} \right.\]
\end{lemma}

\begin{proof}
Let $\mid z \mid =x_2.$
Since
\[ (A_1,A_2,A_4,A_3^{A_4}) \]
has the parameters
$(a_1,a_2,a_4,a_3,x,\tilde{z},y),$
where
\[ z\tilde{z}= \sum_{i=1}^4 a_i^2 +\prod_{i=1}^4a_i+
x^2+y^2-(a_1a_2+a_3a_4)x- (a_1a_4+a_2a_3)y-4\]
we get
\[ z^2\leq  z\tilde{z} \leq (\mid x \mid+4)^2 +(\mid y \mid+4)^2-4.\]
A straight forward computation  yields then the claim.

Similarly  if $\mid y \mid =x_2$ we consider the tuple
\[ (A_1,A_2,A_4^{A_3^{-1}},A_3), \] that has
 the parameters
$(a_1,a_2,a_4,a_3,x,z,\tilde{y}),$
where
\[ y\tilde{y}= \sum_{i=1}^4 a_i^2 +\prod_{i=1}^4a_i+
x^2+z^2-(a_1a_2+a_3a_4)x- (a_1a_3+a_2a_4)z-4.\]
Continuing as above the claim follows.
\end{proof}

\begin{lemma}
For a minimal Fricke tuple holds
 \[ xyz<0  \mbox{ or }  \mid x_2 \mid \leq 12.\]
Moreover, in the case $xyz<0$ the following hold.
\begin{enumerate}
 \item If $x_2\leq x_0+x_1+4$
 then
 $x_2\leq 100$ or $x_0<3$.

 \item
 If $x_0=2$ and  $x_2\leq x_0+x_1+5$ and $a_1,\ldots,a_4 \in \{\pm 2, \pm 1,0\}$ then
 \[ x_2\leq 200\] or
 \[ 0=(a_1-a_2)(a_3-a_4),\; x=2 \]
 or
 \[ 0=(a_1+a_2)(a_3+a_4),\; x=-2. \]
  In particular if $x=2$ the Fricke relation simplifies to
  \[ (y+z-\frac{1}{2}a_1(a_3+a_4))^2+(a_3-a_4)^2(1-\frac{a_1^2}{4})=0. \]
 \end{enumerate}

\end{lemma}

\begin{proof}
Since $\mid a_i \mid \leq 2$ we get
\[ \sum_{i=1}^4 a_i^2 +\prod_{i=1}^4a_i \geq 0.\]
Hence
\[ 4\geq x^2+y^2+z^2+xyz-8(\mid x \mid +\mid y \mid + \mid z \mid )\geq
        (x_2-4)^2-3 \cdot 16+xyz.\]
Therefore the Fricke relation implies
$xyz<0$ or $x_2\leq 12$.
Let $xyz<0$.
\begin{enumerate}
 \item
 Let $x_2> 100.$
 The Fricke relation
\[ \sum_{i=1}^4 a_i^2 +\prod_{i=1}^4a_i+
x^2+y^2+z^2+xyz-(a_1a_2+a_3a_4)x- (a_1a_4+a_2a_3)y-(a_1a_3+a_2a_4)z-4=0,\]
 implies that
 \[ 28+ x^2+y^2+z^2+xyz+8(x_0+x_1+x_2)\geq 0.\]
  Thus
\[ 28+ 2x_0^2+2x_1^2+2x_0x_1+8(x_0+x_1)+16+8(2x_0+2x_1+4) \geq x_0x_1x_2.\]
 Hence
 \[ 76+ 2x_0^2+2x_1^2+2x_0x_1+24(x_0+x_1) \geq x_0x_1^2\]
 and
  \[ 76/x_1^2+6+ 48/x_1 \geq x_0.\]
 Therefore
 $x_0 \leq 11.$
 Thus $x_1\geq x_2-x_0-4 \geq 85.$ Again using
   \[76+ 2x_0^2+2x_1^2+2x_0x_1+24(x_0+x_1) \geq x_0x_1^2\]
 we obtain  $3 > x_0.$
 \item The proof is analogous to (i).
\end{enumerate}

\end{proof}

\begin{rem}\label{comp}
 The restrictions for the minimal Fricke parameters in $\ZZ^7$, 
 where
 $x\neq 2$, allow to search them in the range
 $ [-200\ldots 200]^3 \times [-2 \ldots 2]^4$.
 This can be done with a simple computer program.
\end{rem}

We start with the geometric differential equations that we obtain
from Herfurtner's list. All those are rational pull-backs of
a Gauss hypergeometric differential equation whose monodromy
group is contained in  $\SL_2(\ZZ)$.

In the case where an apparent singularity appears, i.e. the Painlev\'e VI case,
these differential equations were studied by Doran, cf. \cite{dor01}.

\begin{rem}
The following monodromy group generators $\A$ belong to second order
differential equations with $4$ singularities having an apparent singularity.

\[ \begin{array}{llcccc}
Herf.  & (a_1,\ldots,a_4,x,y,z) & A_1&A_2&A_3&A_4 \\

I_1I_1I_2I_2*& (2,2,2,-2,2,0,0) &
\mat{  1}{1 }{ 0}{1  },
&    \mat{ 1}{1 }{ 0}{1  },
&   \mat{ 1}{0 }{ -2}{1  },
&   \mat{ 1}{ -2  }{   2}{ -3 }  \\   
 I_1I_1I_1I_3*&(2,2,2,-2,2,1,-1)&
 \mat{   1}{1 }{ 0}{1  },
&  \mat{  1}{3 }{ 0}{1  },
&   \mat{ 1}{0 }{ -1}{1  },
&   \mat{ 1}{ -4  }{   1}{ -3} \\
I_1I_1I_1III*&(2,2,2,0,1,1,1)&
   \mat{1}{ -1  }{  0}{ 1    },
& \mat{ 1}{ 0  }{  1}{ 1    }, 
& \mat{ 0}{ -1  }{  1}{ 2    },
& \mat{ 1}{ 2  }{  -1}{ -1 }    \\
I_1I_1I_2IV*&(2,2,2,-1,0,0,1)&
  \mat{ 1}{ 1  }{  0}{ 1    },
&\mat{  1}{ 0  }{  -2}{ 1    }, 
& \mat{ 0}{ 1  }{  -1}{ 2 },     
& \mat{ 0}{ -1  }{  1}{ -1}    \\

 I_1I_1 II\; IV*&(2,2,1,-1,1,0,0)&
  \mat{ 1}{ 1  }{  0}{ 1    },
& \mat{  1}{ 0  }{  -1}{ 1    },
& \mat{  0}{ 1  }{  -1}{ 1 },     
& \mat{  0}{ -1  }{ 1}{ -1 }   \\

\end{array} \]

\end{rem}

\begin{proof}
 Since these tuple are in Herfurtner's list they
 arise as pull-backs of hypergeometric differential
 equations. Further their braid group orbit is finite
 and therefore the Fricke parameters are in $[2\ldots 2]^7.$ 
 Using the Herfurtner notation we can easily construct an
 element in the braid group orbit. 
\end{proof}

\begin{tab}
 We list the minimal tuples of monodromy group generators corresponding
 to the Table of geometric Heun equations.
\end{tab}
\[
\begin{array}{ccccccc}
 nr. & (x,y,z,a_1,a_2,a_3,a_4) & A_1 & A_2 & A_3 & A_4 \\
  \mbox{ in }\cite{MR08-2} &&&&&\\
 18&(-5,6, 6,2, 2, -1, -1)&\mat{  1}{ 1  }{  0}{ 1}    ,
                        & \mat{ 1}{ 0  }{  -7}{ 1}    ,
                        & \mat{ 2}{ -1  }{  7}{ -3},      
                        &\mat{ 4}{ -3  }{  7}{ -5}    \\ 

19 &(-10, 5, 5,2, 2, -1, -1)& \mat{  1}{ 2  }{  0}{ 1}    ,
                          & \mat{ 1}{ 0  }{  -6}{ 1}    ,
                          & \mat{ 1}{ -1  }{  3}{ -2},      
                          &\mat{ 4}{ -7  }{  3}{ -5}    \\ 
 20 &(-14,3,11,  2, 2, -1, -1) & \mat{  1}{ 4  }{  0}{ 1}    ,
                         & \mat{ 1}{ 0  }{  -4}{ 1}    , 
                         & \mat{ 1}{ -1  }{  3}{ -2},     
                         & \mat{ 2}{ -7  }{  1}{ -3}     \\
 
 21&  (-4, -4, -3,  2, 2, -1, 0) &\mat{  1}{ 2  }{  0}{ 1}    ,
                           & \mat{ 1}{ 0  }{  -3}{ 1}    , 
                           &\mat{ -1}{ 1  }{  -1}{ 0},      
                           &\mat{ -3}{ 5  }{  -2}{ 3}      \\
  22&(-8, 4, 5,2, 2, -1, 0) & \mat{  1}{ 2  }{  0}{ 1}    ,
                       & \mat{ 1}{ 0  }{  -5}{ 1}    , 
                       &\mat{ 1}{ -1  }{  3}{ -2},      
                       & \mat{ 3}{ -5  }{  2}{ -3}     \\

23 &(-10, 3, 8,2, 2, -1, 0)& \mat{  1}{ 3  }{  0}{ 1}    ,
                       & \mat{ 1}{ 0  }{  -4}{ 1}    , 
                       &\mat{ 1}{ -1  }{  3}{ -2 } ,    
                       &\mat{ 2}{ -5  }{  1}{ -2 }    \\
 
24 & (-3, 4, 6,  2, 2, -1, 1)&   \mat{ 1}{ 1  }{  0}{ 1}    ,
                        & \mat{ 1}{ 0  }{  -5}{ 1}    , 
                        &\mat{ 2}{ -1  }{  7}{ -3} ,     
                        &\mat{ 2}{ -1  }{  3}{ -1}     \\

25 &(-6, 3, 5, 2, 2, -1, 1)& \mat{  1}{ 2  }{  0}{ 1}    ,
                       & \mat{ 1}{ 0  }{  -4}{ 1}    , 
                       &\mat{ 1}{ -1  }{  3}{ -2},     
                       &\mat{ 2}{ -3  }{  1}{ -1}     \\

26 &(-3, -5, -5,2,2,0,0)& \mat{ 1}{ 1  }{  0}{ 1}    ,
                    & \mat{ 1}{ 0  }{  -5}{ 1}    , 
                    & \mat{ -2}{ 1  }{  -5}{ 2},    
                    & \mat{ -3}{ 2  }{  -5}{ 3}     \\

27 &( -6, -4, -4,2,2,0,0) &\mat{  1}{ 2  }{  0}{ 1}    , 
                    &\mat{ 1}{ 0  }{  -4}{ 1}    , 
                    & \mat{ -1}{ 1  }{  -2}{ 1}    , 
                    &\mat{ -3}{ 5  }{  -2}{ 3}     \\

28 &    (-7, -3, -6,2,2,0,0) & \mat{ 1}{3}{ 0}{1 }, 
                     &  \mat{ 1}{ 0  }{  -3}{ 1},    
                     &\mat{ -1}{ 1  }{  -2}{ 1},   
                     & \mat{ -2}{ 5  }{  -1}{ 2}     \\

29 & (-2,-4,-5,2, 2, 0, -1)& 
       \mat{  1}{ 1  }{   0}{ 1 }   ,
     & \mat{    1}{ 0  }{   -4}{ 1}    , 
     & \mat{   -2}{ 1  }{   -5}{ 2},      
     &\mat{-2}{ 1  }{   -3}{ 1 }   

 \\

\end{array}\]
\[
\begin{array}{ccccccc}

30 &(-4,5,6, 2, 2, 0, -1)&   \mat{   1}{ 6  }{   0}{ 1}    , 
    &  \mat{  1}{ 0  }{   -1}{ 1}    ,
    &  \mat{   2}{ -5  }{   1}{ -2},
    &  \mat{   3}{ -13  }{   1}{ -4    }
  \\

31 & (-2,-3,-3,2, 2, -1, -1)&   \mat{    1}{ 2  }{   0}{ 1}    ,  
       &\mat{   1}{ 0  }{   -2}{ 1 }   ,   
       &\mat{  -1}{ 1  }{   -1}{ 0}   ,   
       &\mat{ -2}{ 3  }{  -1}{ 1}   
\\

32 & (3, 2, -2, 2, 1, 0, 0)&\mat{ 1}{ -2  }{  0}{ 1}    ,
                           & \mat{ 0}{ 1  }{  -1}{ 1}    , 
                          &\mat{ 0}{ -1  }{  1}{ 0},      
                          &\mat{ 1}{ 2  }{  -1}{ -1}   \\

 33 & (3, -3, 3,2, 0, 0, 0)& \mat{ 1}{ -3  }{  0}{ 1}    , 
                      &\mat{ 0}{ 1  }{  -1}{ 0}    , 
                      &\mat{ 1}{ 2  }{  -1}{ -1 } ,   
                      &\mat{ -2}{ -5  }{  1}{ 2  }  \\

34 &  (4, 2, -3,2, 1, 0, -1) & \mat{  1}{ -3 }{   0}{ 1},
                & \mat{      0}{ 1 }{   -1}{ 1}, 
               &\mat{    0}{ -1 }{   1}{ 0}, 
               &\mat{     1}{ 3 }{   -1}{ -2}   \\

35 & (-5, -3, -3,2, -1, 1, 1) &\mat{  1}{ -4  }{  0}{ 1}    ,
                         & \mat{ 0}{ -1  }{  1}{ -1}    , 
                         &\mat{ 0}{ -1  }{  1}{ 1   },   
                         &\mat{ -2}{ -7  }{  1}{ 3}    \\

 36 &(3, -4, 4,2, -1, 0, 0)&  \mat{ 1}{ 4  }{  0}{ 1}    ,
                       & \mat{ 0}{ -1  }{  1}{ -1}    , 
                      &\mat{ 2}{ -5  }{  1}{ -2} ,     
                      &\mat{ -3}{ 10  }{  -1}{ 3}    \\

 37& ( -5, 4, 4, 2, 0, -1, -1)&  \mat{ 1}{ 5  }{  0}{ 1}    , 
                         &\mat{ 0}{ 1  }{  -1}{ 0}    , 
                         &\mat{ 1}{ -3  }{  1}{ -2}     
                         &\mat{ 3}{-13  }{  1}{ -4} \\ 

38  & (-5, -5, -5, 2, 1, 1, 1)&  \mat{ 1}{ -6  }{  0}{ 1}    ,
                         & \mat{ 0}{ -1  }{  1}{ 1 }    ,
                         & \mat{ -2}{ -7  }{  1}{ 3 },     
                         &\mat{ -4}{ -21  }{  1}{ 5 }  \\

\end{array}\]

\begin{rem}\label{pullback}
 Let $ (A_1,\ldots,A_4) \in \SL_2(\CC)^4$ with Fricke parameters
 $(x,y,z,0,0,a_3,a_4)$.
then the quadratic pull-back corresponding to
 \[  (A_1,\ldots,A_4) \mapsto (A_3,A_4,A_3^{A_1^{-1}}, A_4^{A_1^{-1}}) \]
 yields the map
\[ (x,y,z,0,0,a_3,a_4) \mapsto (x,-x-yz+a_4\cdot a_3 ,-z^2+2, a_3,a_4,a_3,a_4). \]

  Conversely, if $ (A_1,\ldots,A_4) \in \SL_2(\CC)^4$ with Fricke parameters
 $(x_1,y_1,z_1,a_3,a_4,a_3,a_4)$, where $z_1\neq 2$.
 then $\A$ arises from a tuple
  $ \B \in \SL_2(\CC)^4$ with Fricke parameters
 $(x,y,z,0,0,a_3,a_4)$.
\end{rem}

\begin{lemma}
 We have the following geometric braid group orbits with
 integer Fricke parameters but monodromy group not contained in $\SL_2(\ZZ).$
\[\tiny \begin{array}{cccccc}
  (x,y,z,a_1,a_2,a_3,a_4) & A_1 & A_2 & A_3 & A_4 \\
 (-4, -3, -3, 0, 0, 1, 1) 
    &    \mat{  0}{ -1  }{   1}{ 0    },
&    \mat{  0}{ -2-\sqrt{3}  }{   2-\sqrt{3}}{ 0    },
&    \mat{   1/2+1/2\sqrt{3}}{
          3/2+1/2\sqrt{3}  }{
        -3/2+1/2\sqrt{3}}{
          1/2-1/2\sqrt{3}    },
&     \mat{    1/2+1/2\sqrt{3}}{
          3/2-1/2\sqrt{3}  }{
         -3/2-1/2\sqrt{3}}{
           1/2-1/2\sqrt{3}}    
 \end{array}\]

This group is related to Table~\ref{Her1-18}, row 15, 
via the pull-back corresponding to (cf. Rem.~\label{pullback}
\[  (A_1,\ldots,A_4) \mapsto (A_3,A_4,A_3^{A_1^{-1}},A_4^{A_1^{-1}})\]
\[  (x,y,z,a_1,a_2,a_3,a_4) =(-4, -3, -3, 0, 0, 1, 1) \mapsto (-4,-4,-7,1,1,1,1)\]
and applying the middle convolution
\[  MC_{\zeta_6^{-2}} (\zeta_6 A_3,\zeta_6A_4,\zeta_6 A_3^{A_1},-A_4^{A_1}) \]
 that transforms the Fricke parameters to
\[  (-4,-4,-7,1,1,1,1)\mapsto (-4,-4,-7 ,2,2,2,-1).\]

\end{lemma}

\begin{proof}
Let $\si \in \Gal(\QQ(\sqrt{3}))$.
Then
 \[ (\si(A_1),\ldots,\si(A_4))=(A_1,A_2^{A_1},\si(A_3),\si(A_4)).\]
A monodromy representation in $\SL_2(\ZZ)$ exists if and only if
there exists a $g \in \SL_2( \ZZ(\sqrt{3}))$ such that
 $\si(g)^{-1} g =A_1$.
But this implies that $\si(g_{11}) g_{11}/(\si(g_{12}) g_{12})=-1$.
However there is no element $a=a_1+\sqrt{3}a_2$ in $\QQ(\sqrt{3})$ such that
\[-1=\si(a)a=(a_1+\sqrt{3}a_2)(a_1-\sqrt{3}a_2)= a_1^2-3a_2^2.\]
This can be seen by reducing the equation modulo $3$.

\end{proof}

\begin{lemma} 
 We have the following representatives of
 non geometric braid group orbits with integer Fricke parameters,
 where $\max\{\mid x\mid,\mid y\mid, \mid z\mid\}>2$:
\[ \tiny \begin{array}{llcccc}
 & (a_1,\ldots,a_4,x,y,z) & A_1&A_2&A_3&A_4 \\
 &  (1,1,1,1,2,a^2+2,-a^2-1)&
     \mat{ 0}{ -1  }{   1}{ 1    },
   &  \mat{   1}{ 1  }{   -1}{ 0 },  
   &  \mat{   a}{ -a^2+a-1  }{   1}{ -a+1},      
   &
    \mat{-a+1}{ a^2-a+1  }{   -1}{ a}     \\
 &  (0,0,1,1,2,i(2a-1),-i(2a-1))&
    i\mat{ 1}{ 0  }{   0}{ -1    },
&   -i \mat{ 1}{ 0  }{   0}{ -1    },
& \mat{     a}{ -a^2+a-1  }{   1}{ -a+1 },     
&   \mat{ -a+1}{ a^2-a+1  }{   -1}{ a }    \\
 &  (0,0,0,0,2,2a,-2a)&
   i \mat{ 1}{ 0  }{   0}{ -1    },
&  -i\mat{   1}{ 0  }{   0}{ -1    },
 &  i  \mat{ a}{ 1  }{   -a^2-1}{ -a    },
&  -i \mat{  -a}{ -a^2-1  }{  1}{ a}    \\
&(2,2,0,0,2,a,-a)
& \mat{ 1}{ -a  }{   0}{ 1    },
& \mat{ 1}{ a  }{   0}{ 1    },
& \mat{ 0}{ 1  }{   -1}{0    },
& \mat{ 0}{ -1  }{   1}{ 0    } \\
&(2,2,1,1,2,a+1,-a+1)
& \mat{ 1}{ -a  }{   0}{ 1    },
& \mat{ 1}{ a  }{   0}{ 1    },
& \mat{ 0}{ 1  }{   -1}{1    },
& \mat{ 1}{ 1  }{   -1}{ 0    } \\
&(2,2,0,1,2,-a,a+1)
& \mat{ 1}{ a  }{   0}{ 1    },
& \mat{ 1}{ -a-1  }{   0}{ 1    },
& \mat{ 0}{ 1  }{   -1}{0    },
& \mat{ 0}{ -1  }{   1}{ 1    } \\

\end{array}\]
\end{lemma}

\begin{proof}
 Let $(x,y,z,a_1,\ldots,a_4)=(x,y,z,1,1,1,1).$
 Then 
 \[    MC_{\zeta_3}(-\zeta_3 A_1,-\zeta_3 A_2,-\zeta_3 A_3,-A_4)=(B_1,\ldots,B_4) \in \SL_2(\CC)^4 \] 
 where by Lemma~\ref{prop} and Cor.~\ref{MCFricke}
 \[  (x,y,z,\tr(B_1),\ldots, \tr(B_4))=(x,y,z,2,2,2,-1).\]
 Since the middle convolution preserves geometric differential equations
 the claim for the first case follows from Cor.~\ref{Nongeom2}.
 The second case can be reduced to the first case via a pull-back
 as in Rem.~\ref{pullback}.
 The third case is settled by Theorem~\ref{Chud}. 
 The fourth case can be reduced to the case $(a_1,\ldots,a_4)=(2,2,2,2)$
 via a quadratic  pull-back as in Rem.~\ref{pullback} and the claim follows from
 Cor.~\ref{Nongeom2}.   
 In the fifth case we consider the braid group action
 \[ (A_1,\ldots,A_4) \mapsto (A_1,A_2^{A_3},A_3^{A_2A_3},A_4).\]
 Since $\tr(A_1A_2^{A_3})=a^2+a+2>2$ for $a\not \in \{-1,0\}$ the application
 $MC_{-1}$ to $\A$ yields a tuple $\B$ in $\GO_4(\ZZ)^4$ by Lemma~\ref{prop}.
 By taking a quadratic pull-back that changes $\B$ to $(B_3,B_4,B_3^{B_1},B_4^{-1})$
 the claim follows
 from Corollary~\ref{nongeom} ii) and Lemma~\ref{prop} ii).
 Again as in the previous case the claim follows for the last case
 from Corollary~\ref{nongeom} and the  application of
 $MC_{-1}$.
\end{proof}

\begin{rem}\label{isgeom}
 The only remaining cases are those where we can't decide
 whether a  braid group orbit in $\SL_2(\ZZ)^4$ is
 geometric or not.
 (We list the minimal tuples)
\[
\begin{array}{ccccccc}
  & (x,y,z,a_1,a_2,a_3,a_4) & A_1 & A_2 & A_3 & A_4 \\
  & (-1, 4, 4,2, 2, 1, 1) &  \mat{    1}{ 3  }{   0}{ 1 }   , 
      &    \mat{1}{ 0  }{   -1}{ 1 }   , 
      &    \mat{ 2}{ -3  }{   1}{ -1 },     
   & \mat{2}{ -3  }{   1}{ -1  }   \\

  & (3, -1, 3, 2, 1, 1, 1) &     \mat{  1}{ 2  }{   0}{ 1 }   ,
                           &    \mat{ 0}{ -1  }{   1}{ 1 }   , 
                           &  \mat{  0}{ -1  }{   1}{ 1 },     
  & \mat{ 0}{ 1  }{   -1}{ 1 }    \\

  &( 2, -1, 2,2, 1, 1, 0)& 
                          \mat{1}{ 1  }{   0}{ 1}    ,
                         & \mat{    0}{ -1  }{   1}{ 1}    , 
                        &  \mat{  0}{ -1  }{   1}{ 1}, 
                        &  \mat{  0}{1  }{  -1}{ 0 }  
\end{array} \]
\end{rem}

{
\begin{tab}\label{geomHeun}
 Geometric Heun equations whose monodromy group is after scaling in $\SL_2(\ZZ)$ (cf.\cite{MR08-2}):
\[y''+( \frac{1-\th_1}{x-t}+\frac{1-\th_2}{x}+\frac{1-\th_3}{x-1})y'+\frac{\th_{41}\th_{42}x-q}{x(x-1)(x-t)}y=0\]
 \small
\begin{center}
\begin{tabular}{|c|c|c|c|c|c|c|c|c|}
\hline
*  &   $q$ & $t$ & $\th_1$ & $\th_2$ & $\th_3$   & $\th_{42}$ & $\th_{41}$  \\ \hline
\hline
1 &
$ 
\frac{-1}{3}t_1 $ &
$\frac{t_1^2}{3}, \ t_1^2+3t_1+3=0 $ &
$  0  $ &
$  0  $ &
$  0  $ &
$  1  $ &
$  1  $

\\ \hline
2 &
$  0  $ &
$ -1 $ &
$  0  $ &
$  0  $ &
$  0  $ &
$  1 $ &
$  1  $ \\ \hline
3 &
$  -2 $ &
$ -8 $ &
$  0 $ &
$  0  $ &
$  0  $ &
$  1  $ &
$  1  $ \\  \hline
4 &
$  -3t_1 $ &
$-t_1^2,\ t_1^2-11 t_1-1=0$ &
$  0 $ &
$  0  $ &
$  0  $ &
$  1  $ &
$  1 $
\\ \hline
5 &
$  0  $ &
$ -1 $ &
$  0  $ &
$  0  $ &
$  0  $ &
$  1  $ &
$  1 $
\\ \hline
6 &
$
\frac{-1}{3}t_1$ &
$\frac{t_1^2}{3}, t_1^2+3t_1+3=0$ &
$  0  $ &
$  0 $ &
$  0 $ &
$  1  $ &
$  1  $
\\ \hline
7 & 
$   \frac{-23}{243}t_1 $ &
$\frac{t_1^2}{9},3t_1^2-14t_1+27=0$ &
$ 0 $ &
$ \frac{1}{3} $ &
$  0  $ &
$  \frac{5}{6}  $ &
$  \frac{5}{6}  $
\\ \hline
8 &
$ 
\frac{731}{1152} $ &
$ \frac{81}{32} $ &
$ 0  $ &
$ \frac{1}{3} $ &
$ 0  $ &
$  \frac{5}{6}  $ &
$  \frac{5}{6}  $
\\ \hline
9 &
$
\frac{-125}{12}  $ &
$ -80 $ &
$  0  $ &
$  0  $ &
$ \frac{1}{3} $ &
$ \frac{5}{6} $ &
$ \frac{5}{6}  $
\\ \hline
10 &
$ 
\frac{-25}{18}$ &
$ -\frac{27}{5} $ &
$ \frac{1}{3} $ &
$  0  $ &
$  0  $ &
$ \frac{5}{6} $ &
$ \frac{5}{6}  $
\\ \hline
11 &
$ 
\frac{-25}{512}t_1 $ &
$\frac{t_1^2}{8}, 4t_1^2+13t_1+32=0$ &
$  0  $ &
$ \frac{1}{2} $ &
$  0  $ &
$  \frac{3}{4}  $ &
$  \frac{3}{4}  $
\\ \hline
12 & 

$
\frac{9}{64}$ &
$ \frac{1}{4} $ &
$  0  $ &
$ \frac{1}{2} $ &
$ 0   $ &
$ \frac{3}{4} $ &
$ \frac{3}{4}  $
\\ \hline
13 &
$
 \frac{39}{500}  $ &
$ -\frac{3}{125} $ &
$ \frac{1}{2} $ &
$ 0  $ &
$ 0  $ &
$  \frac{3}{4}  $ &
$  \frac{3}{4}  $
\\ \hline
14 &
$ 
\frac{-3}{4}$ &
$ -3 $ &
$ \frac{1}{2} $ &
$  0 $ &
$  0  $ &
$  \frac{3}{4}    $ &
$ \frac{3}{4}  $
\\ \hline
15 &
$  0  $ &
$ -1 $ &
$  0  $ &
$ \frac{2}{3} $ &
$  0  $ &
$  \frac{2}{3}  $ &
$  \frac{2}{3}   $
\\ \hline
16 &
$  \frac{7}{9}  $ &
$ \frac{27}{2} $ &
$  0  $ &
$ \frac{2}{3} $ &
$  0  $ &
$  \frac{2}{3}  $ &
$  \frac{2}{3}   $
\\ \hline
17 & 
$
\frac{-2}{9} $ &
$ -1 $ &
$ \frac{2}{3} $ &
$  0  $ &
$  0  $ &
$  \frac{2}{3}   $ &
$  \frac{2}{3} $
\\ \hline

18 &
$
\frac{1}{21}t_1 $ &
$\frac{t_1^2}{49},\ t_1^2-13 t_1+49=0$ &
$ \frac{1}{3} $ &
$ 0  $ &
$ \frac{1}{3} $ &
$   \frac{2}{3}  $ &
$  \frac{2}{3}  $
\\ \hline
19 &

$ 0  $ &
$ -1 $ &
$\frac{ 1}{3} $ &
$   0  $ &
$ \frac{1}{3} $ &
$  \frac{2}{3}  $ &
$  \frac{2}{3}  $
\\ \hline
20 & 
$\frac{-2}{3}t_1 $ &
$-\frac{t_1^2}{2}, t_1^2-10t_1-2$ &
$  0  $ &
$ \frac{1}{3} $ &
$  0  $ &
$   1  $ &
$  \frac{2}{3}  $
\\ \hline
21 & 

$ 
\frac{1}{144} (78\zeta+43)
$ & 
$ -\frac{2}{7}  (3 \zeta+1),\ \zeta^2+3=0 $ &
$ 0  $ &
$ \frac{1}{2} $ &
$\frac{1}{3} $ &
$  \frac{7}{12}  $ &
$ \frac{7}{12}$ 
\\ \hline
22 &

$ 
\frac{721}{2250}  $ &
$ \frac{189}{125} $ &
$\frac{ 1}{2} $ &
$ \frac{1}{3} $ &
$  0  $ &
$ \frac{7}{12} $ &
$ \frac{7}{12} $
\\ \hline
23 &

$ 
\frac{77}{972} $ &
$ -\frac{1}{27 } $ &
$ \frac{1}{2} $ &
$ 0 $ &
$ \frac{1}{3} $ &
$ \frac{7}{12} $ &
$ \frac{7}{12} $
\\ \hline
24 &
$  -\frac{1}{36}  $ &
$ -\frac{16}{9} $ &
$  0  $ &
$ \frac{2}{3} $ &
$ \frac{1}{3} $ &
$ \frac{1}{2} $ &
$ \frac{1}{2} $
\\ \hline
25 &
$  -\frac{1}{2} $ &
$ 9 $ &
$ \frac{1}{3} $ &
$\frac{ 2}{3} $ &
$  0  $ &
$ \frac{1}{2} $ &
$ \frac{1}{2} $
\\ \hline
26 &
$
\frac{1}{25}t_1 $ &
$\frac{4t_1^2}{125},\ t_1^2-11t_1+125/4=0$ &
$ \frac{1}{2} $ &
$  0  $ &
$ \frac{1}{2} $ &
$  \frac{1}{2}  $ &
$  \frac{1}{2} $
\\ \hline
27 & 
$  0 $ &
$ -1 $ &
$ \frac{1}{2} $ &
$  0  $ &
$ \frac{1}{2} $ &
$ \frac{1}{2} $ &
$ \frac{1}{2}$
\\ \hline
28 &
$
\frac{-1}{6}t_1$ &
$-\frac{t_1^2}{3},\ t_1^2-6t_1-3=0$ &
$  0  $ &
$ \frac{1}{2} $ &
$  0  $ &
$  1  $ &
$  \frac{1}{2}  $
\\ \hline
29 &
$  \frac{-5}{972} $ &
$ -\frac{5}{27} $ &
$ \frac{1}{2} $ &
$ \frac{2}{3} $ &
$  0  $ &
$  \frac{5}{12}  $ &
$  \frac{5}{12}  $
\\ \hline
30 & 
$   \frac{5}{18}$ &
$ 5 $ &
$ \frac{1}{2} $ &
$ \frac{2}{3} $ &
$  0  $ &
$  \frac{5}{12}  $ &
$  \frac{5}{12}  $
\\ \hline
31 &
$  0 $ &
$ -1 $ &
$\frac{ 2}{3} $ &
$  0  $ &
$ \frac{2}{3} $ &
$ \frac{1}{3} $ &
$ \frac{1}{3}  $
\\ \hline
32 & 
$  0  $ &
$ -1 $ &
$ \frac{1}{2} $ &
$ \frac{2}{3} $ &
$ \frac{1}{2} $ &
$  \frac{1}{6}  $ &
$  \frac{1}{6}  $
\\ \hline
33 &
$ 
\frac{-1}{48}t_1$ &
$\frac{t_1^2}{3},\ t_1^2+3t_1+3=0$ &
$ \frac{1}{2}  $ &
$\frac{1}{2}   $ &
$\frac{1}{2}  $ &
$\frac{1}{4}  $ &
$  \frac{1}{4}  $
\\ \hline
34 &
$  0 $ &
$ -\frac{1}{3} $ &
$ \frac{1}{2} $ &
$ \frac{2}{3} $ &
$ \frac{1}{3} $ &
$  \frac{1}{4}  $ &
$  \frac{1}{4}  $
\\ \hline
35 &

$  0  $ &
$ -1 $ &
$ \frac{1}{3} $ &
$\frac{ 2}{3} $ &
$ \frac{1}{3} $ &
$  \frac{1}{3}  $ &
$  \frac{1}{3}  $
\\ \hline 
36 &
$
\frac{4}{243}t_1 $ &
$\frac{t_1^2}{27},\ t_1^2-10t_1+27=0$ &
$ \frac{1}{2} $ &
$ \frac{1}{3} $ &
$ \frac{1}{2 }$ &
$  \frac{1}{3}  $ &
$  \frac{1}{3}  $
\\ \hline
37 &
$ 
\frac{-25}{3072}t_1 $ &
$\frac{t_1^2}{64},t_1^2+11 t_1+64=0$ &
$ \frac{1}{3} $ &
$\frac{ 1}{2} $ &
$ \frac{1}{3} $ &
$  \frac{5}{12}  $ &
$  \frac{5}{12}  $
\\ \hline
38 &
$ 
\frac{-1}{12}t_1  $ &
$\frac{t_1^2}{3},\ t_1^2+3t_1+3=0$ &
$ \frac{1}{3} $ &
$\frac{ 1}{3} $ &
$\frac{ 1}{3} $ &
$ \frac{ 1}{2}   $ &
$  \frac{ 1}{2} $
\\ \hline
\end{tabular}
\end{center}
\end{tab}
}

{}

\begin{center}
\vspace{.25in} 
{{\sc Author's adress:
 Stefan Reiter}} \\          ̈
Universit\"at Mainz, Fachbereich 17, Mathematik, 55099 Mainz, Germany\\
E-mail:
{\tt reiters@uni-mainz.de} 
\end{center}

\end{document}